\documentclass[12pt]{article}

\usepackage{amssymb}

\newtheorem{thm}{Theorem}[section]
\newtheorem{defn}[thm]{Definition}
\newtheorem{prop}[thm]{Proposition}
\newtheorem{cor}[thm]{Corollary}
\newtheorem{lemma}[thm]{Lemma}
\newtheorem{rema}[thm]{Remark}

\newcommand{\halmos}{\rule{1ex}{1.4ex}}

\newcommand{\bea}{\begin{eqnarray}}
\newcommand{\eea}{\end{eqnarray}}
\newcommand{\nn}{\nonumber \\}
\newcommand{\be}{\begin {equation}}
\newcommand{\pa}{\partial}
\newcommand{\ee}{\end{equation}}

\newcommand{\cc}{\varphi}

 \newcommand{\res}{\mbox{\rm Res}}
\renewcommand{\hom}{\mbox{\rm Hom}}
 \newcommand{\pf}{{\it Proof.}\hspace{2ex}}
 \newcommand{\epf}{\hspace*{\fill}\mbox{$\halmos$}}
 \newcommand{\epfv}{\hspace*{\fill}\mbox{$\halmos$}\vspace{1em}}
 \newcommand{\epfe}{\hspace{2em}\halmos}

\newcommand{\wt}{\mbox{\rm wt}\ }

\newcommand{\lbar}{\bigg\vert}

\title{ {\bf Intertwining operator superalgebras 
and vertex tensor categories
for superconformal algebras, I} }
\date{}
\author{Yi-Zhi Huang and Antun Milas}

\begin{document}

\bibliographystyle{alpha}
\maketitle

\begin{abstract}
We  construct the intertwining operator superalgebras 
and vertex tensor
categories for the $N=1$ superconformal minimal models and 
other related
models. 
\end{abstract}

\renewcommand{\theequation}{\thesection.\arabic{equation}}
\renewcommand{\thethm}{\thesection.\arabic{thm}}
\setcounter{equation}{0}
\setcounter{thm}{0}
\setcounter{section}{-1}

\section{Introduction}

Superconformal symmetries play a fundamental role in string 
theory since perturbative string theory
is described by two-dimensional quantum field theories with 
superconformal symmetries.
On the other hand, superconformal algebras also arise 
in mathematical structures, for example,
in the quantum field theory structure underlying
the moonshine module for the Monster finite simple group (see
\cite{FLM1},
\cite{B}, \cite{FLM2}, \cite{DGM}, \cite{DGH} and \cite{H2.1}) and 
in the study of mirror symmetry for 
Calabi-Yau manifolds (see \cite{G1}, \cite{G2}, \cite{D}, 
\cite{LVW},
\cite{GP} and many other works by physicists discussed in the survey
\cite{Gr}). A complete mathematical understanding of superconformal
field theories is needed in order to 
solve the related mathematical problems.

Intertwining operator (super)algebras and vertex tensor
categories are  equivalent essentially 
to genus-zero weakly-holomorphic
conformal field theories in the sense of Segal \cite{S1} \cite{S2} satisfying
additional properties
(see \cite{H2.5} and \cite{H4}). 
The general theory of
intertwining operator algebras, including a construction of 
intertwining operator algebras {}from representations of suitable 
vertex operator algebras, was developed by the first author 
in \cite{H2.5}, \cite{H4} and \cite{H6}. 
A tensor product theory for modules 
for a vertex operator algebra, including a construction of 
vertex tensor categories {}from categories of modules for 
vertex operator algebras, was
developed by Lepowsky and the first author in \cite{HL1}--\cite{HL6}
and \cite{H1}.
Besides giving  genus-zero weakly-holomorphic
conformal field theories, the theory of 
intertwining operator algebras is also closely related to the
so called ``boundary conformal field theories'' first developed by Cardy
in \cite{C1}
and \cite{C2}, which have many applications in condensed matter physics 
and have recently been  applied successfully 
(see, for example, \cite{RS} and \cite{FS}) to the study of 
D-branes in string theory.
We expect that the theory of intertwining operator algebras
will provide a mathematical foundation to a number of ``world sheet''
constructions in 
conformal field theories, boundary conformal field theories and 
D-branes. 

In the present series of papers, we
shall construct intertwining operator superalgebras and vertex tensor
categories associated to the superconformal minimal models 
and other related
models.

In this paper (Part I), we apply previously obtained results on
representations of $N=1$ superconformal algebras and associated vertex
operator algebras to show that the general theory for the construction
of intertwining operator (super)algebras and the tensor product theory
for modules for a vertex operator (super)algebra are applicable in
this case. Therefore we obtain an intertwining operator superalgebra
structure on the direct sum of all inequivalent irreducible modules
for a minimal $N=1$ superconformal vertex operator superalgebra. We
also obtain a vertex tensor category structure and consequently a
braided tensor category structure associated to these models.  The
main work in this paper is to prove that the conditions to use the
theory of intertwining operator algebras and the tensor product theory
are satisfied for these models. These results are also generalized
easily to a much more general class of vertex operator superalgebras.

The present paper is organized as follows: In Section 1, we recall 
the notion of $N=1$ superconformal vertex operator superalgebra.
In Section 2, we recall and prove
some basic results on representations
of minimal $N=1$ superconformal vertex operator superalgebras
and of vertex operator superalgebras in a much more general class.
 Section 3 is devoted to 
the proof of the convergence and extension properties for 
products of intertwining operators for 
minimal $N=1$ superconformal vertex operator superalgebras and 
for vertex operator superalgebras in the general class. 
The main tool is null vectors and differential equations.
Our main results on
the intertwining operator superalgebra structures and 
vertex tensor category structures are given in Section 4. 
In the appendix (Section 5), 
for one particular example, we give explicitly the calculations of 
the fusion rules, null vectors and  differential 
equations.

{\bf Acknowledgment}: The research of 
Y.-Z. H. is supported in part by
NSF grant DMS-9622961.

\renewcommand{\theequation}{\thesection.\arabic{equation}}
\renewcommand{\thethm}{\thesection.\arabic{thm}}
\setcounter{equation}{0}
\setcounter{thm}{0}

\section {$N=1$ superconformal vertex operator superalgebras}

In this section we recall the notion of $N=1$ superconformal vertex
operator algebra and basic properties of such an algebra. These
algebras have been studied extensively by physicists, and
are formulated precisely by
Kac and Wang \cite{KW} and by Barron \cite{B1} \cite{B2}.

\begin{defn}
{\rm An {\it $N=1$ superconformal vertex operator superalgebra} is a
vertex operator superalgebra $(V, Y, \mathbf{1}, \omega)$ together with
an odd element $\tau$ called the {\it Neveu-Schwarz element} 
satisfying the
following axiom: Let
$$Y(\tau, x)=\sum_{n\in \mathbb{Z}}G(n+1/2)x^{-n-2}.$$
Then   the following 
$N=1$ Neveu-Schwarz relations hold: For $m, n\in \mathbb{Z}$,
\begin{eqnarray*}
{[L(m), L(n)]}&=&(m-n)L(m+n)+\frac{c}{12}(m^{3}-m)\delta_{m+n, 0},\\
{[L(m), G(n+ 1/2)]}&=&\left(\frac{m}{2}-(n+
1/2)\right)G(m+n+ 1/2),\\
{[G(m+1/2), G(n-1/2)]}&=& 
2L(m+n)+\frac{c}{3}(m^{2}+m)\delta_{m+n, 0},
\end{eqnarray*}
where $L(m)$, $m\in \mathbb{Z}$, are the Virasoro operators on $V$ and
$c$ is the central charge of $V$.

{\it Modules} and {\it intertwining operators} for an $N=1$ superconformal 
vertex operator superalgebra are modules and intertwining operators
for the underlying vertex operator superalgebra.}
\end{defn}

The $N=1$ superconformal vertex operator superalgebra defined above is
denoted by $(V, Y, \mathbf{1}, \tau)$ (without $\omega$ since 
$\omega=L(-2)\mathbf{1}=\frac{1}{2}G(-1/2)\tau$) 
or simply $V$. 
Note that
a module $W$ for a vertex operator superalgebra 
(in particular the algebra itself) 
has a
$\mathbb{Z}_{2}$-grading called {\it sign} in addition to the
$\mathbb{C}$-grading by weights. We shall always use 
$W^{0}$ and $W^{1}$ to denote the even and odd subspaces of 
$W$. If $W$ is irreducible, there exists
$h\in \mathbb{C}$ such that $W=W^{0}\oplus W^{1}$ where
$W^{0}=\coprod_{n\in h+\mathbb{Z}}W_{(n)}$ and $W^{1}=\coprod_{n\in
h+\mathbb{Z}+1/2}W_{(n)}$ are the even and odd parts of $W$,
respectively. We shall always use the notation $|\cdot|$ to denote the
map from the union of the even and odd subspaces of a
vertex operator superalgebra or of 
a module for such an algebra to $\mathbb{Z}_{2}$ 
by taking the signs of elements in the union.

The notion of $N=1$ superconformal vertex operator superalgebra above
was reformulated using odd formal variables by Barron.  For a
complete and detailed discussion, see \cite{B1} and \cite{B2}. 
Here we only give the
parts we need in later sections.

For $l$ symbols
$\varphi_{1}, \dots, \varphi_{l}$, consider
the exterior algebra of
the vector space over $\mathbb{C}$ spanned by these symbols.  
We denote this exterior algebra by $\mathbb{C}[\cc_{1}, \dots,
\cc_{l}]$. For any 
vector space $E$, we use $E[\cc_{1}, \dots,
\cc_{l}]$ to denote the tensor product of $E$ and $\mathbb{C}[\cc_{1}, \dots,
\cc_{l}]$. In particular, if $E$ is the polynomial algebra 
$\mathbb{C}[x_{1}, \dots, x_{k}]$ generated by formal variables
$x_{1}, \dots, x_{k}$ or the space
$\mathbb{C}[[x_{1}, \dots, x_{k}]]$ of formal Laurent series 
generated by these formal variables, we have the 
vector space $\mathbb{C}[x_{1}, \dots, x_{k}][\cc_{1}, \dots,
\cc_{l}]$ or $\mathbb{C}[[x_{1}, \dots, x_{k}]][\cc_{1}, \dots,
\cc_{l}]$. In this case, we call $x_{1}, \dots, x_{k}$ and $\cc_{1},
\dots, \cc_{l}$ {\it even} and {\it odd} formal variables,
respectively.   Note that any element of $\mathbb{C}[x_{1}, \dots, x_{k}]
[\cc_{1}, \dots, \cc_{l}]$ is a linear combination of monomials in
$x_{1}, \dots, x_{k}$ and $\varphi_{1}, \dots, \varphi_{l}$.  For each
monomial, the total order in $x_{1}, \dots, x_{k}$ and 
$\varphi_{1}, \dots, \varphi_{l}$
and the order in $\varphi_{1}, \dots, \varphi_{l}$
modulo $2$ give a $\mathbb{Z}$-grading called {\it degree} and a
$\mathbb{Z}_{2}$-grading called {\it sign}, respectively,
to $\mathbb{C}[x_{1}, \dots,
x_{k}][\cc_{1}, \dots, \cc_{l}]$.  With the $\mathbb{Z}_{2}$-grading,
$\mathbb{C}[x_{1}, \dots, x_{k}][\cc_{1}, \dots, \cc_{l}]$ is an
associative superalgebra. 
Similarly,  $\mathbb{C}[[x_{1}, \dots, x_{k}]][\cc_{1}, \dots,
\cc_{l}]$
is also a superalgebra. 

For any vector space $E$, consider the vector space
\begin{eqnarray*}
&E[x_{1}, \dots, x_{k}][\cc_{1}, \dots, \cc_{l}],&\\
&E[x_{1}, x_{1}^{-1}, 
\dots, x_{k}, x_{k}^{-1}][\cc_{1}, \dots, \cc_{l}],&\\
&E[[x_{1}, \dots, x_{k}]][\cc_{1}, \dots, \cc_{l}],&\\
&E[[x_{1}, x_{1}^{-1}, 
\dots, x_{k}, x_{k}^{-1}]][\cc_{1}, \dots, \cc_{l}],&\\
&E\{x_{1}, 
\dots, x_{k}\}[\cc_{1}, \dots, \cc_{l}]&
\end{eqnarray*}
and
$$E((x_{1}, 
\dots, x_{k}))[\cc_{1}, \dots, \cc_{l}].$$
If $E$ is a $\mathbb{Z}_{2}$-graded vector space, then there 
are  natural structures of modules over the ring 
$\mathbb{C}[x_{1}, \dots, x_{k}][\cc_{1}, \dots, \cc_{l}]$
on these spaces.

Let $(V, Y, \mathbf{1}, \tau)$ be an 
$N=1$ superconformal vertex operator superalgebra. 
We define the {\it vertex operator map with odd variable}
\begin{eqnarray*}
Y: V\otimes V&\to &V((x))[\varphi]\\
u\otimes v&\mapsto& Y(u, (x, \varphi)) v
\end{eqnarray*}
by
$$Y(u, (x, \varphi))v=Y(u, x)v+\varphi Y(G(-1/2)u, x)v$$
for $u, v\in V$. (We use the same notation $Y$ to denote 
the vertex operator map and the vertex operator map with odd variable.)
Then we have:

\begin{prop}\label{odd-svoa}
The vertex operator map with odd variable satisfies the following 
properties:

\begin{enumerate}

\item The vacuum property: 
$$Y(\mathbf{1}, (x, \cc))=1$$
where $1$ on the right-hand side is the identity map on $V$.

\item The creation property: For any $v\in V$,
$$ Y(v, (x,\cc))\mathbf{1} \in V[[x]][\cc], $$
$$ \lim_{(x, \cc) \mapsto (0, 0)} Y(v, (x,\cc))\mathbf{1}=v.$$

\item The Jacobi identity: In $(\mbox{\rm End}\ V)[[x_{0}, x_{0}^{-1}, 
x_{1}, x_{1}^{-1}, 
x_{2}, x_{2}^{-1}]][\cc_{1}, \cc_{2}]]$, we have 
\bea
\lefteqn{x_0^{-1} \delta \left 
(\frac {x_1-x_2-\cc_{1}\cc_{2}}{x_0} \right)
Y(u, (x_1,\cc_{1})) Y(v, (x_2,\cc_{2}))} \nn
&&\quad -(-1)^{|u||v|}x_0^{-1} \delta \left ( \frac
{x_2-x_1+
\cc_{1}
\cc_{2}}{-x_0} \right ) Y(v, (x_2,\cc_{2})) Y(u, (x_1,\cc_{1}))
\nn
&&= x_2^{-1} \delta \left ( \frac {x_1-x_0- \cc_{1}
\cc_{2}}{x_2} \right
)Y(Y(u, (x_0,\cc_{1}-\cc_{2}))v, (x_2,\cc_{2}))
\eea
for $u, v\in V$ which are either even or odd.

\item The $G(-1/2)$-derivative property: For any $v\in V$,
$$Y(G(-1/2) v, (x, \cc))=\left(\frac{\partial}{\partial \cc}+
\cc\frac{\partial}{\partial x}\right)Y( v, (x,\cc)),$$

\item The $L(-1)$-derivative property: For any $v\in V$, 
$$Y(L(-1) v, (x,\cc))=\frac{\partial}{\partial x}
Y(v, (x,\cc)).$$

\item The {\it skew-symmetry}: For any $u, v\in V$ which are either 
even or odd,
\bea 
Y(u,
(x,\varphi))v=(-1)^{|u||v|}e^{xL(-1)+\varphi G(-1/2)}Y(v,
(-x,-\varphi))u.  \epfe
\eea

\end{enumerate}

\end{prop}

The proof of this result is straightforward  and can be found in \cite{B2}.

We can also reformulate the data and axioms for modules 
and intertwining operators for an $N=1$ 
superconformal vertex operator superalgebra using odd variables. 
Since the goal of the present paper is to construct an algebra {}from
intertwining operators for certain $N=1$ 
superconformal vertex operator superalgebras, we give the details of 
the corresponding reformulation of the data and axioms 
for intertwining operators using odd variables. 

Let $W_{1}$, $W_{2}$ and $W_{3}$ be modules for an $N=1$ 
superconformal vertex operator superalgebra $V$ and $\mathcal{Y}$ an 
intertwining operator of type ${W_{3}\choose W_{1}W_{2}}$. We define the 
corresponding {\it intertwining operator map with odd variable}
\begin{eqnarray*}
\mathcal{Y}: W_{1}\otimes W_{2}&\to &W_{3}\{x\}[\varphi]\\
w_{(1)}\otimes w_{(2)}&\mapsto& \mathcal{Y}(w_{(1)}, (x, \varphi)) w_{(2)}
\end{eqnarray*}
by
$$\mathcal{Y}(w_{(1)}, (x, \varphi)) w_{(2)}
=\mathcal{Y}(w_{(1)}, x) w_{(2)}+
\varphi\mathcal{Y}(G(-1/2)w_{(1)}, x) w_{(2)}$$
for $u, v\in V$.
Then we have:

\begin{prop}
The intertwining operator map with odd variable satisfies the following 
properties:

\begin{enumerate}

\item The {\it Jacobi identity}: In $\hom(W_{1}\otimes W_{2}, W_{3})
\{x_{0},
x_{1},
x_{2}\}[\cc_{1}, \cc_{2}]$, we have 
\bea
\lefteqn{x_0^{-1} \delta \left ( \frac {x_1-x_2- \cc_{1}\cc_{2}}{x_0} \right ) 
Y(u, (x_1,\cc_{1})) \mathcal{Y}(w_{(1)}, (x_2,\cc_{2}))} \nn
&& -(-1)^{|u||w_{(1)}|}x_0^{-1} \delta \left ( \frac
{x_2-x_1+
\cc_{1}
\cc_{2}}{-x_0} \right ) \mathcal{Y}(w_{(1)}, (x_2,\cc_{2})) 
Y(u, (x_1,\cc_{1}))
\nn
&&= x_2^{-1} \delta \left ( \frac {x_1-x_0- \cc_{1}
\cc_{2}}{x_2} \right
)\mathcal{Y}(Y(u, (x_0,\cc_{1}-\cc_{2}))w_{(1)}, (x_2,\cc_{2}))
\eea
for $u\in V$ and $w_{(1)}\in W_{1}$ which are either even or odd.

\item The {\it $G(-1/2)$-derivative property}: For any $v\in V$,
$$\mathcal{Y}(G(-1/2)
 w_{(1)}, (x, \cc))=\left(\frac{\partial}{\partial \cc}+
\cc\frac{\partial}{\partial x}\right)
\mathcal{Y}(w_{(1)}, (x,\cc)),$$

\item The {\it $L(-1)$-derivative property}: For any $v\in V$, 
$$\mathcal{Y}(L(-1) w_{(1)}, (x,\cc))=\frac{\partial}{\partial x}
\mathcal{Y}(w_{(1)}, (x,\cc)).$$

\item The {\it skew-symmetry}: There is a linear
isomorphism
$$\Omega: \mathcal{V}^{W_{3}}_{W_{1} W_{2}} \to
\mathcal{V}^{W_{3}}_{W_{2} W_{1}}
$$
such that 
\begin{eqnarray*}
\lefteqn{\Omega (\mathcal{ Y})(w_{(1)}, (x,\cc))w_{(2)}}\nn
&&=(-1)^{|w_{(1)}||w_{(2)}|}e^{xL(-1)+\cc 
G(-1/2) }\mathcal{ Y}(w_{(2)}, (e^{-\pi i}x,-\cc)))w_{(1)}
\end{eqnarray*}
for $w_{(1)}\in W_{1}$ and $w_{(2)}\in W_{2}$ which are either even or odd. \epf

\end{enumerate}

\end{prop}

The proof of this result is similar to the proof of Proposition 
\ref{odd-svoa} and is omitted.

\renewcommand{\theequation}{\thesection.\arabic{equation}}
\renewcommand{\thethm}{\thesection.\arabic{thm}}
\setcounter{equation}{0}
\setcounter{thm}{0}

\section{Minimal $N=1$ superconformal vertex operator superalgebras}

In this section, we recall the constructions and results on minimal
$N=1$ superconformal vertex operator superalgebras and their
representations. Some new results needed in later
sections are also proved.  
We then introduce in this section a class of vertex
operator superalgebras and generalize most of the results for minimal
$N=1$ superconformal vertex operator superalgebras to algebras in this
class. A large part of the material in this section is {}from
\cite{KW} and \cite{A}.

The {\it $N=1$ Neveu-Schwarz Lie superalgebra} 
is the Lie superalgebra 
$$\mathfrak{n}\mathfrak{s}^{(1)}=\oplus_{n\in \mathbb{Z}}\mathbb{C}L_{n}\oplus
\oplus_{n\in \mathbb{Z}}G_{n+1/2}\oplus \mathbb{C}C$$
satisfying the following 
$N=1$ Neveu-Schwarz relations:
\begin{eqnarray*}
{[L_{m}, L_{n}]}&=&(m-n)L_{m+n}+\frac{C}{12}
(m^{3}-m)\delta_{m+n, 0},\\
{[L_{m}, G_{n+ 1/2}]}&=&\left(\frac{m}{2}-\left(n+
\frac{1}{2}\right)\right)G_{m+n+ 1/2},\\
{[G_{m+1/2}, G_{n-1/2}]}&=& 
2L_{m+n}+\frac{C}{3}(m^{2}+m)\delta_{m+n, 0},\\
{[C, L_{m}]}&=&0,\\
{[C, G_{m+1/2}]}&=&0
\end{eqnarray*}
for $m, n\in \mathbb{Z}$. For simplicity, we shall simply denote the $N=1$
Neveu-Schwarz Lie superalgebra by $\mathfrak{n}\mathfrak{s}$ in this paper.

Note that the elements $L_0,L_{1}, L_{-1}, G_{1/2}, G_{-1/2}$ of
$\mathfrak{n}\mathfrak{s}$ span a subalgebra.  It is known that 
this subalgebra is isomorphic to
$\mathfrak{o}\mathfrak{s}\mathfrak{p}(2,1)$.

We now construct representations of the $N=1$ Neveu-Schwarz 
Lie superalgebra.
Consider the two subalgebras 
\begin{eqnarray*} 
\mathfrak{n}\mathfrak{s}^{+}&=&\oplus_{n>0}\mathbb{C}L_{n}\oplus
\oplus_{n\ge 0}G_{n+1/2},\\
\mathfrak{n}\mathfrak{s}^{-}&=&\oplus_{n<0}\mathbb{C}L_{n}\oplus
\oplus_{n<0}G_{n+1/2}
\end{eqnarray*}
of $\mathfrak{n}\mathfrak{s}$. 
Let $U(\cdot)$ be the functor {}from the 
category of Lie superalgebras to
the category of associative algebras obtained by taking the universal 
enveloping algebras of  Lie superalgebras. 
For any representation of $\mathfrak{n}\mathfrak{s}$, 
we shall use 
$L(m)$ and $G(m+1/2)$, $m\in \mathbb{Z}$,
to denote the representation images of $L_{m}$ and 
$G_{m+1/2}$.
For any $c, h\in \mathbb{C}$, the Verma module $M_{\mathfrak{ns}}(c, h)$ 
for $\mathfrak{n}\mathfrak{s}$ 
is a free $U(\mathfrak{n}\mathfrak{s}^{-})$-module generated by
$\mathbf{1}_{c, h}$ such that
\begin{eqnarray*}
\mathfrak{n}\mathfrak{s}^{+}\mathbf{1}_{c, h}&=&0,\\
L(0)\mathbf{1}_{c, h}&=&h\mathbf{1}_{c, h},\\
C\mathbf{1}_{c, h}&=&c\mathbf{1}_{c, h}.
\end{eqnarray*}
There exists a unique maximal proper submodule $J_{\mathfrak{ns}}(c,
h)$
of $M_{\mathfrak{ns}}(c, h)$. It 
is easy to see that when $c\ne 0$, 
$\mathbf{1}_{c, 0}$, $G(-3/2)\mathbf{1}_{c, 0}$ and
$L(-2)\mathbf{1}_{c, 0}$ are not in $J_{\mathfrak{n}\mathfrak{s}}(c, 0)$.
Let 
$$
L_{\mathfrak{ns}}(c, h)=M_{\mathfrak{ns}}(c, h)/J_{\mathfrak{ns}}(c, h)
$$
and 
$$V_{\mathfrak{ns}}(c, 0)=M_{\mathfrak{ns}}(c, 0)/\langle 
G(-1/2)\mathbf{1}_{c, 0}\rangle$$
where $\langle 
G(-1/2)\mathbf{1}_{c, 0}\rangle$ is the submodule of 
$M_{\mathfrak{ns}}(c, 0)$ generated by $G(-1/2)\mathbf{1}_{c, 0}$.
Then $L_{\mathfrak{ns}}(c, 0)$ and $V_{\mathfrak{ns}}(c, 0)$
have the structures of vertex operator superalgebras with
the vacuum $\mathbf{1}_{c, 0}$, the Neveu-Schwarz element
$G(-3/2)\mathbf{1}_{c, 0}$ and the Virasoro element
$L(-2)\mathbf{1}_{c, 0}$ (see
\cite{KW}).

The following result was conjectured by Kac and Wang 
in \cite{KW} and proved by Adamovi\'{c} in \cite{A}
using the relationship between representations of 
the $N=1$ Neveu-Schwarz Lie superalgebra and representations of the affine 
Lie algebra $A_{1}^{(1)}$ on the rational level 
obtained in \cite{GKO}, \cite{KWak1}, 
\cite{KWak2}, \cite{AM} and \cite{As}:

\begin{thm}\label{ad}
The vertex operator superalgebra $L_{\mathfrak{ns}}(c, 0)$ 
has finitely many irreducible 
modules and every module for $L_{\mathfrak{ns}}(c, 0)$ 
is completely reducible if
and only if 
$$
c=c_{p, q}=\frac{3}{2}\left(1-2\frac{(p-q)^{2}}{pq}\right)
$$
where $p, q$ are integers larger than $1$ such that 
$p-q\in 2\mathbb{Z}$ and $(p-q)/2$ and $q$ are 
 relatively prime to each other.
A set of representatives of the equivalence classes of irreducible modules
for $L_{\mathfrak{ns}}(c_{p, q}, 0)$ is 
$$
\{L_{\mathfrak{ns}}(c_{p, q}, h_{p, q}^{m, n})\}_{0<m<p, 
\; 0<n<q, \; m, n\in \mathbb{Z},
\; m-n\in 2\mathbb{Z}}
$$
where for any $m, n\in \mathbb{Z}$ satisfying $0<m<p, 0<n<q$ and
$m-n\in 2\mathbb{Z}$,
$$
h_{p, q}^{m, n}=\frac{(np-mq)^{2}-(p-q)^{2}}{8pq}. \hspace{\fill}
\epfe
$$
\end{thm}

For any pair $p, q$ of integers larger than $1$ such that 
$p-q\in 2\mathbb{Z}$ and $(p-q)/2$ and $q$ are 
 relatively prime to each other, 
we call the vertex operator algebra $L_{\mathfrak{ns}}(c_{p, q}, 0)$ a
{\it minimal $N=1$ superconformal vertex operator superalgebra}.

\begin{prop}\label{fusion}
Let $m_{i}, n_{i}\in \mathbb{Z}$, $i=1, 2, 3$, satisfying $0<m_{i}<p,
0<n_{i}<q$ and $m_{i}-n_{i}\in 2\mathbb{Z}$ and
$\mathcal{ Y}$ an intertwining operator of type
$${L_{\mathfrak{n}\mathfrak{s}}(c_{p, q}, h_{p,
q}^{m_{3}, n_{3}})\choose L_{\mathfrak{n}\mathfrak{s}}(c_{p, q}, h_{p,
q}^{m_{1}, n_{1}})L_{\mathfrak{n}\mathfrak{s}}(c_{p, q}, h_{p,
q}^{m_{2}, n_{2}})}.$$  Then we have:

\begin{enumerate}
\item\label{2.2-1}  For any $w_{(1)}\in L_{\mathfrak{n}\mathfrak{s}}(c_{p, q}, h_{p,
q}^{m_{1}, n_{1}})$ and $w_{(2)}\in L_{\mathfrak{n}\mathfrak{s}}(c_{p, q}, h_{p,
q}^{m_{2}, n_{2}})$, 
$$\mathcal{Y}(w_{(1)}, x)w_{(2)}\in x^{h_{p,
q}^{m_{3}, n_{3}}-h_{p,
q}^{m_{1}, n_{1}}-h_{p,
q}^{m_{2}, n_{2}}}L_{\mathfrak{n}\mathfrak{s}}(c_{p, q}, h_{p,
q}^{m_{3}, n_{3}})((x^{1/2}))$$

\item\label{2.2-2} The map $\mathcal{Y}$ is uniquely 
determined by the maps
$$(\mathbf{1}_{c, h_{p,
q}^{m_{1}, n_{1}}})_{h_{p,
q}^{m_{1}, n_{1}}+h_{p,
q}^{m_{2}, n_{2}}-h_{p,
q}^{m_{3}, n_{3}}-1}$$ and 
$$(G(-1/2)\mathbf{1}_{c, h_{p,
q}^{m_{1}, n_{1}}})_{h_{p,
q}^{m_{1}, n_{1}}+h_{p,
q}^{m_{2}, n_{2}}-h_{p,
q}^{m_{3}, n_{3}}-1/2}$$
from $W_{2}$ to $W_{3}$
(recalling  that $\mathbf{1}_{c, h_{p,
q}^{m_{1}, n_{1}}}$ is the lowest weight
vector in $L_{\mathfrak{n}\mathfrak{s}}(c_{p, q}, h_{p,
q}^{m_{1}, n_{1}})$), that is,
\bea
\lefteqn{(\mathbf{1}_{c, h_{p,
q}^{m_{1}, n_{1}}})_{h_{p,
q}^{m_{1}, n_{1}}+h_{p,
q}^{m_{2}, n_{2}}-h_{p,
q}^{m_{3}, n_{3}}-1}}\nn
&&=(G(-1/2)\mathbf{1}_{c, h_{p,
q}^{m_{1}, n_{1}}})_{h_{p,
q}^{m_{1}, n_{1}}+h_{p,
q}^{m_{2}, n_{2}}-h_{p,
q}^{m_{3}, n_{3}}-1/2}\nn
&&=0
\eea
implies
$\mathcal{Y}=0$.

\item\label{2.2-3} The space 
$$\mathcal{V}^{L_{\mathfrak{n}\mathfrak{s}}(c_{p, q}, h_{p,
q}^{m_{3}, n_{3}})}_{L_{\mathfrak{n}\mathfrak{s}}(c_{p, q}, h_{p,
q}^{m_{1}, n_{1}})L_{\mathfrak{n}\mathfrak{s}}(c_{p, q}, h_{p,
q}^{m_{2}, n_{2}})}$$
is at most 2-dimensional.

\end{enumerate}
\end{prop}
\pf
Conclusion \ref{2.2-1} is clear since the three modules are 
irreducible.

We prove Conclusion \ref{2.2-2} now. For convenience, we denote 
$h_{p,
q}^{m_{3}, n_{3}}-h_{p,
q}^{m_{1}, n_{1}}-h_{p,
q}^{m_{2}, n_{2}}$ 
by $\Delta$ and $\mathbf{1}_{c, h_{p,q}^{m_{i}, n_{i}}}$ 
for $i=1, 2, 3$ by $w_{(i)}$, respectively. Suppose that  
$$(w_{(1)})_{\Delta-1}=(G(-1/2)w_{(1)})_{\Delta-1/2}=0$$
but $\mathcal{Y}\neq 0$.

We need the following commutator formulas which are 
consequences of the Jacobi identity: 
For $n \in \mathbb{Z}, m \in \Delta+\mathbb{Z}/2$,
\bea
\lefteqn{[L(n),(w_{(1)})_m]}\nn
&&=(-m-n-1+(n+1)h_{p,
q}^{m_{1}, n_{1}})(w_{(1)})_{m+n},\\
\lefteqn{[L(n), (G(-1/2)w_{(1)})_{m}]}\nn
&&=(-m-n-1+(n+1)h_{p, q}^{m_{1}, n_{1}})(G(-1/2)w_{(1)})_{m+n},\\
\lefteqn{[G(n+1/2),(w_{(1)})_m]}\nn
&&=(G(-1/2)w_{(1)})_{m+n+1},\\
\lefteqn{[G(n+1/2),(G(-1/2)w_{(1)})_{m}]}\nn
&&=(-m-n-1/2-2h_{p,q}^{m_{1},n_{1}})(w_{(1)})_{m+n}.
\eea

We claim  that either 
$$\mathcal{Y}(w_{(1)}, x)w_{(2)} \neq 0$$
or 
$$\mathcal{Y}(G(-1/2)w_{(1)}, x)w_{(2)} \neq 0.$$ 
In fact, 
if 
$$\mathcal{Y}(w_{(1)}, x)w_{(2)}=\mathcal{Y}(G(-1/2)w_{(1)}, x)w_{(2)}=0,$$
then by the commutator formulas above,
$\mathcal{Y}(w_{(1)}, x)=0$ which is a contradiction
since $L_{\mathfrak{n}\mathfrak{s}}(c_{p, q}, h_{p,
q}^{m_{1}, n_{1}})$ is an irreducible 
$\mathfrak{n}\mathfrak{s}$-module 
and $\mathcal{Y}\ne 0$.

Let
$$k=\max \{m\in \Delta +\mathbb{Z}/2\;|\; 
(w_{(1)})_m w_{(2)} \neq 0 \ \mbox{or} 
\ (G(-1/2)w_{(1)})_{m+1/2}w_{(2)} \neq 0\}.$$
If  
$$(w_{(1)})_k w_{(2)} \neq 0,$$
then
{}from the commutator formulas above, it follows that
$$L(n)(w_{(1)})_k w_{(2)}=(-k-n-1+(n+1)h_{p, q}^{m_{1}, n_{1}})
(w_{(1)})_{k+n}w_{(2)}=0$$
for $n \geq 1$ and 
$$G(n+1/2)(w_{(1)})_k
w_{(2)}=(G(-1/2)w_{(1)})_{k+n+1} w_{(2)}=0$$
 for $n \geq 0$. Thus 
$(w_{(1)})_k w_{(2)}$ is a multiple of the lowest weight vector 
$w_{(3)}$.
Then
$$h_{p,
q}^{m_{3}, n_{3}}=\wt ((w_{(1)})_k w_{(2)})=h_{p,
q}^{m_{1}, n_{1}}-k-1+h_{p,
q}^{m_{2}, n_{2}}.$$ 
It follows that $k=\Delta-1$. So  $(w_{(1)})_{k}w_{(2)}=0$, a contradiction.

If 
$$(G(-1/2)w_{(1)})_{k+1/2} w_{(2)} \neq 0,$$ 
then 
\begin{eqnarray*}
\lefteqn{L(n)(G(-1/2)w_{(1)})_{k+1/2} w_{(2)}}\nn
&&=(-k-1/2-n-1+(n+1)h_{p, q}^{m_{1}, n_{1}})(G(-1/2)w_{(1)})_{k+1/2+n}\nn
&&=0
\end{eqnarray*}
for $n\ge 1$ and 
\begin{eqnarray*}
\lefteqn{G(n+1/2)(G(-1/2)w_{(1)})_{k+1/2} w_{(2)}}\nn
&&=(-k-n-1-2h_{p, q}^{m_{1}, n_{1}}(w_{(1)})_{k+n+1}\nn
&&=0,\end{eqnarray*}
for
$n \geq 0$. Thus 
$(G(-1/2)w_{(1)})_{k+1/2} w_{(2)}$ is a highest
weight vector and therefore 
$$h_{p,
q}^{m_{3}, n_{3}}=\wt (G(-1/2)w_{(1)})_{k+1/2} w_{(2)}=h_{p,
q}^{m_{1}, n_{1}}-k-1+h_{p,
q}^{m_{2}, n_{2}}.$$ 
It follows that $k=\Delta-1$. So $(G(-1/2)w_{(1)})_{k+1/2} w_{(2)}=0$, 
a contradiction.

Conclusion \ref{2.2-3} follows immediately {}from 
Conclusion \ref{2.2-2}.  
\epfv

Combining Theorem \ref{ad} and the third conclusion of 
Proposition \ref{fusion}, we obtain:

\begin{cor}
The minimal 
$N=1$ superconformal vertex operator superalgebras are rational in the 
sense of \cite{HL1}, that is, the following three conditions are 
satisfied: 

\begin{enumerate}

\item Every module for such an algebra is completely 
reducible.

\item There are only finitely many inequivalent irreducible modules for 
such an algebra.

\item The fusion rules among any three (irreducible) modules are 
finite.\epf

\end{enumerate}
\end{cor}

\begin{rema}
{\rm In \cite{Z}, Zhu introduced a weaker notion of module for a
vertex operator algebra and a notion of rational vertex operator
algebra. Zhu's notion of rational vertex operator algebra is different
{}from the notion of rational vertex operator algebra in \cite{HL1}
because Zhu's notion requires a stronger completely reducibility
result for modules in his sense. 
The complete reducibility result and the classification of 
irreducible modules proved in \cite{A} together with the 
third conclusion of Proposition \ref{fusion}
gives only the rationality in the sense of \cite{HL1}.}
\end{rema}

We also have:

\begin{prop}
Any finitely-generated 
lower truncated generalized $L_{\mathfrak{n}\mathfrak{s}}(c_{p,q},0)$-module
$W$ is an ordinary module.
\end{prop}
\pf
Suppose that $W$ is 
generated by a single vector $w \in W$. Then
by the Poincar\'{e}-Birkhoff-Witt theorem and 
the lower truncation condition, every
homogeneous  subspace of $U(\mathfrak{n}\mathfrak{s})w$ is 
finite-dimensional, proving the result.
\epfv

Let $n$ be a positive integer, $(p_{i}, q_{i})$, $i=1, \dots, n$, $n$
pairs of integers larger than $1$ such that $p_{i}-q_{i}\in
2\mathbb{Z}$ and $(p_{i}-q_{i})/2$ and $q_{i}$ are relatively prime to
each other, and let 
$V=L_{\mathfrak{ns}}(c_{p_{1}, q_{1}}, 0)\otimes \cdots \otimes
L_{\mathfrak{ns}}(c_{p_{n}, q_{n}}, 0)$. 
{}From the trivial generalizations of the
results proved in \cite{FHL} and \cite{DMZ} to vertex operator
superalgebras, $V$ is a rational $N=1$ superconformal 
vertex operator superalgebra, a set
of representatives of equivalence classes of irreducible modules for
$V$ can be listed 
explicitly and the fusion rules for $V$ are finite and 
can be calculated easily.

We introduce a class of $N=1$ superconformal vertex operator vertex
operator superalgebras:

\begin{defn}
{\rm Let $n$ be a positive integer, $(p_{i}, q_{i})$, 
$i=1, \dots, n$, $n$
pairs of integers larger than $1$ such that $p_{i}-q_{i}\in
2\mathbb{Z}$ and $(p_{i}-q_{i})/2$ and $q_{i}$ are relatively prime to
each other.  An $N=1$ superconformal vertex operator vertex
operator superalgebra $V$ is
said to be {\it in the class $\mathcal{C}_{p_{1}, q_{1};
\dots; p_{n}, q_{n}}$}
if $V$ has a vertex operator subalgebra 
isomorphic to $L_{\mathfrak{n}\mathfrak{s}}(c_{p_{1}, q_{1}},
0)\otimes \cdots \otimes L_{\mathfrak{n}\mathfrak{s}}(c_{p_{m}, q_{m}}, 0)$.}
\end{defn}

\begin{prop}\label{2-7}
Let $V$ be an $N=1$ superconformal vertex operator vertex
operator superalgebra in the class $\mathcal{C}_{p_{1}, q_{1};
\dots; p_{n}, q_{n}}$. Then any finitely-generated 
lower truncated generalized $V$-module
$W$ is an ordinary module.
\end{prop}
\pf
The proof is similar to the proof of Proposition 3.7 in \cite{H2}.
So here we only point out the main difference. As in \cite{H2}, 
we discuss only the case $n=2$. Similar to 
the proof of Proposition 3.7 in \cite{H2}, using the Jacobi identity,
the formula $G(-1/2)^{2}=L(-1)$ and Theorem 4.7.4 of \cite{FHL},
 we can reduce our 
proof in the case of $n=2$ 
to the finite-dimensionality of the space spanned 
by the elements of the form 
\begin{eqnarray}\label{generators}
\lefteqn{L(-m^{(1)}_{1})\cdots L(-m^{(1)}_{e_{1}})
G(-a^{(1)}_{1})\cdots G(-a^{(1)}_{r_{1}})
(L(-1)^{l_{1}}G(-1/2)^{k_{1}}
u_{(j)}^{(1)})_{j_{1}}\cdot}\nn
&&\quad\cdot
L(n^{(1)}_{1})\cdots L(n^{(1)}_{f_{1}})
G(b^{(1)}_{1})\cdots G(b^{(1)}_{s_{1}})w_{(t)}^{1}\nn
&&\quad\otimes L(-m^{(2)}_{1})
\cdots L(-m^{(2)}_{e_{2}})
G(-a^{(2)}_{1})\cdots G(-a^{(2)}_{r_{2}})\cdot\nn
&&\quad\cdot (L(-1)^{l_{2}}G(-1/2)^{k_{2}}u_{(j)}^{(2)})_{j_{2}}
L(n^{(2)}_{1})
\cdots L(n^{(2)}_{f_{2}})G(b^{(2)}_{1})\cdots G(b^{(2)}_{s_{2}})
w_{(t)}^{2},\nn
&&
\end{eqnarray}
for $m^{(1)}_{1}, \dots, m^{(1)}_{e_{1}}$, 
$m^{(2)}_{1}, \dots, m^{(2)}_{e_{2}}$,
$n^{(1)}_{1}, \dots, n^{(1)}_{f_{1}}$, 
$n^{(2)}_{1}, \dots, n^{(2)}_{f_{2}}$, 
$a^{(1)}_{1}, \dots, a^{(1)}_{r_{1}}$, 
$a^{(2)}_{1}, \dots, a^{(2)}_{r_{2}}$,
$b^{(1)}_{1}, \dots, b^{(1)}_{s_{1}}$, 
$b^{(2)}_{1}, \dots, b^{(2)}_{s_{2}}
\in \mathbb{Z}_{+}$,
$l_{1}, l_{2}\in \mathbb{N}$, $k_{1}, k_{2}=0, 1$,
$j_{1}, j_{2}\in \mathbb{Q}$, $t=1, \dots, c$,
$j=1, \dots, d$, where $u^{(i)}_{(j)}$, $j=1, \dots, d$, $i=1, 2$,
are elements of $V$ such that the 
$L_{\mathfrak{n}{\mathfrak{s}}}(c_{p_{i}, q_{i}}, 0)$-submodules 
generated by them isomorphic to 
$L(c_{p, q}, h_{p_{i}, q_{i}}^{m_{j}, n_{j}})$ with the images of 
$u^{(i)}_{(j)}$, $j=1, \dots, d$, $i=1, 2$, as the lowest weight vectors 
and such that $V$ is isomorphic to the direct sum of these submodules,
and where
$w_{(t)}^{(i)}$, $t=1, \dots, c$, $i=1, 2$, are homogeneous
elements of some irreducible 
$L_{\mathfrak{n}\mathfrak{s}}(c_{p_{i}, q_{i}}, 0)$-modules. 
Using the $L(-1)$-derivative property for generalized modules,
we see that elements of the form (\ref{generators}) are spanned by 
elements of the form (\ref{generators}) 
with $l_{1}=l_{2}=0$. There are only 
finitely many of elements of the form (\ref{generators}) 
with $l_{1}=l_{2}=0$
and of a fixed weight $s$ because $W$ is lower truncated. 
Thus the homogeneous subspaces of 
$W$ are finite dimensional. So $W$ is a $V$-module. 
\epfv

\renewcommand{\theequation}{\thesection.\arabic{equation}}
\renewcommand{\thethm}{\thesection.\arabic{thm}}
\setcounter{equation}{0}
\setcounter{thm}{0}

\section{Correlation functions and differential equations of regular
singular points}

In this section, we study products and iterates of intertwining
operators for the minimal $N=1$ superconformal vertex operator
superalgebras. The goal is to prove that these products and iterates
satisfying the convergence and extension properties introduced in
\cite{H1}. The main tool is differential equations 
of regular singular
points. 
Though the strategy of the proof is similar to that in \cite{H2}, 
there are 
subtle and nontrivial differences. We shall be brief on the parts of 
proofs which are similar to those in \cite{H2} but give detailed discussions
on the parts which are different.

Let $p, q$ be integers larger than $1$ such that 
$p-q\in 2\mathbb{Z}$ and $(p-q)/2$ and $q$ are 
 relatively prime to each other. Let $W_{i}$, 
$i=1, \dots, 5$, be irreducible modules for the vertex operator algebra 
$L_{\mathfrak{n}\mathfrak{s}}(c_{p, q}, 0)$ and 
$\mathcal{Y}_1$ and $\mathcal{Y}_2$ 
intertwining operators of type ${W_{4}\choose{W_{1} W_{5}}}$
and ${W_{5}\choose{W_{2} W_{3}}}$, respectively.

We first state a proposition which 
describes matrix coefficients of the 
products of $\mathcal{Y}_1$ and $\mathcal{Y}_2$ with
odd variables. The proof is very easy and is omitted.

\begin{prop} \label{struc}
Let $w_{(i)} \in W_i$ (i=1,2,3) and $w'_{(4)}\in W'_{4}$. Then
\bea
\lefteqn{\langle w'_{(4)}, \mathcal{Y}_{1}(w_{(1)}, 
(x_{1},\varphi_1))
\mathcal{Y}_{2}(w_{(2)}, (x_{2},\varphi_2))w_{(3)} \rangle} \nn
&&= Q_{0,0}(x_1,x_2)+\varphi_1 Q_{0,1}(x_1,x_2)
+\varphi_2Q_{1,0}(x_1,x_2)+
\varphi_1\varphi_2 Q_{1,1}(x_1,x_2)\nn
\eea
where $Q_{k,l}(x_1,x_2)=x_1^{t_1-k/2}x_2^{t_2-l/2}
R_{k,l}(x_2^{1/2}/x_1^{1/2})$,
and $t_1,t_2 \in \mathbb{Q}$, $R_{k,l}(x) \in \mathbb{C}[[x]]$ 
($k, l=0,1$).\epf
\end{prop}

We first prove the following theorem:

\begin{thm} \label{top}
Let $w_{(i)} \in W_i$ (i=1,2,3) and $w'_{(4)}\in W'_{4}$
be the lowest weight vectors. Then we have:

\begin{enumerate}

\item The series 
$$ 
\langle w'_{(4)}, \mathcal{Y}_{1}(w_{(1)}, (x_{1},\varphi_1))
\mathcal{Y}_{2}(w_{(2)}, (x_{2},\varphi_2))w_{(3)} \rangle
|_{x_j^n=e^{n\log z_j}, j=1, 2} 
$$
is convergent to 
a multivalued (Grassman) analytic function in the region  $|z_1|>|z_2|>0$ for any choice of $\log z_{1}$ and $\log z_{2}$. 

\item For each $k, l=0,1$  there
exist $m_{k,l} \in
\mathbb{N}$, analytic functions $f_{i}^{k,l}(z)$ in the region 
$|z|<1$ and  $s_{i}^{k, l}, r_{i}^{k, l} \in \mathbb{Q}$ for 
$i=1, \dots, m_{k, l}$,  such that 
$Q_{k,l}(x_1,x_2)|_{x_j^n=e^{n\log z_j}, j=1, 2}$
can be analytically extended to a multivalued analytic functions of the form 
\bea \label{ext}
\sum_{i=1}^{m_{k,l}} z_2^{s_{i}^{k,l}}(z_1-z_2)^{r_{i}^{k,l}}
f^{k,l}_i\left(\frac{z_1-z_2}{z_2}\right),
\eea
when $|z_2|>|z_1-z_2|>0$.
\end{enumerate}
\end{thm}
%
%
\pf
{}From Proposition \ref{struc}, we see that to prove the
Conclusion 1, it is enough to prove that the series
$Q_{k,l}(x_1.x_2)|_{x_i^n=e^{n\log z_i}, i=1, 2}$, $k, l=0, 1$,
 are absolutely convergent 
 in the region $|z_1|>|z_2|>0$.

We need  the following structural result about singular vectors
obtained by Astashkevich in \cite{As}:

\begin{lemma} \label{Asing}
Let $M_{\mathfrak{n}\mathfrak{s}}(c_{p.q},h)$ be a Verma module for $\mathfrak{n}\mathfrak{s}$
and $w \in M_{\mathfrak{n}\mathfrak{s}}(c_{p.q},h)$ a singular vector of weight $n+h$ ($n \in \mathbb{N}/2$). Then 
\bea
w=a(G(-1/2)^{2n}+\cdots)\mathbf{1}_{c_{p, q}, h},
\eea
where $a$ is a nonzero complex number and 
$\cdots$ denotes a sum of monomials different {}from $G(-1/2)^{2n}$. 
\epf
\end{lemma}
 
Since the irreducible module $W_{3}$ is
a quotient of a Verma module 
by the submodule generated by two singular vectors,
we can find nonzero $P\in U(\mathfrak{n}\mathfrak{s}^{-})$ such that 
$Pw_{(3)}=0$ in $W_{3}$. Moreover, we can choose such a 
$P$ to be even 
because otherwise we can compose it with $G(-1/2)$. 
By Lemma \ref{Asing}, we can alway normalize $P$ such that
$P=G(-1/2)^{2n}+\cdots$ 
where  
$\cdots$ denotes a sum of monomials different {}from $G(-1/2)^{2n}$.

Now we prove Conclusion 1 
for $Q_{0,1}$ and $Q_{1,0}$.
{}From $Pw_{(3)}=0$ and the commutator formula
between $G(n+1/2)$, 
$L(n)$, $n\in \mathbb{Z}$, and 
intertwining operators, we obtain a system 
of differential equations
\bea 
0&=&\langle w'_{(4)}, \mathcal{ Y}_{1}(G(-1/2)w_{(1)}, z_{1})
\mathcal{ Y}_{2}(w_{(2)}, z_{2})Pw_{(3)}
\rangle\nn
&=&D_1(z_1,z_2)Q_{1,0}+D_2(z_1,z_2)Q_{0,1},\label{system1} \\
0&=&\langle w'_{(4)}, \mathcal{ Y}_{1}(w_{(1)}, z_{1})
\mathcal{ Y}_{2}(G(-1/2)w_{(2)}, z_{2})Pw_{(3)}
\rangle\nn
&=&D_3(z_1,z_2)Q_{1,0}+D_4(z_1,z_2)Q_{0,1}\label{system1-1}
\eea
for $Q_{1,0}$ and $Q_{0,1}$ of regular 
singular points with the only possible 
singular points $z_1, z_{2}=0, \infty$. 

Since
$P=G(-1/2)^{2n}+\cdots$ where $\cdots$ denotes a sum of  
monomials different 
{}from $G(-1/2)^{2n}$, we have
\begin{eqnarray*}
D_1&=&(\partial_{z_1}+\partial_{z_2})^n+ \sum_{i+j 
\leq n-1}f_{i,j}(z_1,z_2)\partial^i_{z_1}\partial^j_{z_2},\nn
D_2&=&\sum_{i+j \leq n-1}g_{i,j}(z_1,z_2)\partial^i_{z_1}\partial^j_{z_2},\nn
D_3&=&\sum_{i+j
\leq n-1}k_{i,j}(z_1,z_2)\partial^i_{z_1}\partial^j_{z_2},\nn
D_4&=&(\partial_{z_1}+\partial_{z_2})^n+ \sum_{i+j 
\leq n-1}h_{i,j}(z_1,z_2)\partial^i_{z_1}\partial^j_{z_2}
\end{eqnarray*}
for certain meromorphic
functions $f_{i,j},g_{i,j},h_{i,j}$ and $k_{i,j}$.

Because of the structure of $Q_{1,0}$ and $Q_{0,1}$, we can
reduce the system (\ref{system1})--(\ref{system1-1}) to a system 
\begin{eqnarray*}
\hat{D}_1(z)R_{1,0}(z)+\hat{D}_2(z)R_{0,1}(z)&=&0
 \\
\hat{D}_3(z)R_{1,0}(z)+\hat{D}_4(z)R_{0,1}(z)&=&0
\end{eqnarray*}
with
analytic coefficients in one variable $z=z_2/z_{1}$ in the region
$0<|z|<1$ (there might be a singularity at $z=0$) for $R_{1,0}$ and
$R_{0,1}$,
where $\hat{D}_i(z)$, $i=1, 2, 3, 4$, have the following form 
\begin{eqnarray*}
\hat{D}_1(z)&=&\partial_z^n+ \sum_{i \leq n-1} f_i(z) \partial_z^i,
\\
\hat{D}_2(z)&=&\sum_{i \leq n-1} g_i(z) \partial_z^i,
\\
\hat{D}_4(z)&=&\partial_z^n+ \sum_{i \leq n-1} h_i(z) \partial_z^i, 
\\
\hat{D}_3(z)&=&\sum_{i \leq n-1} k_i(z) \partial_z^i
\end{eqnarray*}
where $f_i,g_i,h_i,k_i$ are analytic functions with possible poles at $z=0$.
This (in general higher order) linear system of equations can be further reduced
to a linear system of equations of regular singular points
of  first order.
(More precisely we introduce new unknowns $A_i=R_{1,0}^{(i)}$ 
and $B_i=R_{0,1}^{(i)}$,
$i=0, \dots, n-1$. Then we have a first-order linear 
system of $2n$ equations for $A_i$'s and
$B_i$'s, $i=0, \dots, n-1$.)
{}From the theory of differential equations
of regular singular points,
it follows that $R_{1,0}(z)$ and $R_{0,1}(z)$ are
absolutely convergent in the region $|z|<1$.  
Thus $Q_{1, 0}$ and $Q_{0,1}$ are absolutely convergent to 
analytic
functions in the region $|z_{1}|>|z_{2}|>0$. 

Now we would like to extend $Q_{1,0}$ and $Q_{0, 1}$ analytically
to functions of the form (\ref{ext}) in the region 
$|z_{2}|>|z_{1}-z_{2}|>0$.
Let $Q$ be a nonzero
element of $U(\mathfrak{n}\mathfrak{s}^{-})$ 
such that $Qw_{(2)}=0$ in $W_{2}$. Similar to $P$, we can choose 
$Q$ to be even and we can  normalize
$Q$ such that $Q=G(-1/2)^{2m}+\cdots$ where $\cdots$ denotes a sum of  
monomials different 
{}from $G(-1/2)^{2m}$.
{}From $Qw_{(2)}=0$, the commutator formulas between 
$G(n+1/2)$ and $L(n)$, $n\in \mathbb{Z}$, and
the Jacobi identity for intertwining operators, we have
\bea 
0&=&\langle w'_{(4)}, \mathcal{ Y}_{1}(G(-1/2)w_{(1)}, z_{1})
\mathcal{ Y}_{2}(Qw_{(2)}, z_{2})w_{(3)}
\rangle\nn
&=&D'_1(z_1,z_2)Q_{1,0}+D'_2(z_1,z_2)Q_{0,1},
\label{system4} \\
0&=&\langle w'_{(4)}, \mathcal{ Y}_{1}(w_{(1)}, z_{1})
\mathcal{ Y}_{2}(G(-1/2)Qw_{(2)}, z_{2})Pw_{(3)}
\rangle\nn
&=&D'_3(z_1,z_2)Q'_{1,0}+D_4(z_1,z_2)Q_{0,1}\label{system4-1}
\eea
for some differential operators $D'_{1}$, $D'_{2}$, 
$D'_3$ and $D'_4$.

{}From the formulas used to derive 
(\ref{system4})--(\ref{system4-1}) and from 
$Q=G(-1/2)^{2m}+\cdots$ where $\cdots$ denotes a sum of  
monomials different 
{}from $G(-1/2)^{2m}$, 
it is easy to see that  (\ref{system4})--(\ref{system4-1})
is a system of equations of 
regular singular points with the only possible
singular points $z_{2}, z_{1}-z_{2}=0, \infty$. 
Since (\ref{system1})--(\ref{system1-1}) is a system of 
equations of regular singular points with the only possible
singular points $z_{1}=z_{2}=0, \infty$, the equations
(\ref{system1}), (\ref{system1-1}) and the equations
(\ref{system4}), (\ref{system4-1}) are independent.
Thus we obtain a system (\ref{system1})--(\ref{system4-1})  of equations 
of regular singular points with the only possible singularities
$z_{1}, z_{2}=0, \infty$ and $z_{1}=z_{2}$. 
This system is consistent because it has a formal series solution 
$(Q_{1, 0}, Q_{0,1})$. 
 Since we have proved that 
$Q_{1,0}$ and $Q_{0,1}$ are absolutely convergent in the region 
$|z_{1}|>|z_{2}|>0$ and since the intersection 
$|z_{1}|>|z_{2}|>|z_{1}-z_{2}|>0$
of the regions 
$|z_{1}|>|z_{2}|>0$ and $|z_{2}|>|z_{1}-z_{2}|>0$ are nonempty, 
we can use the sum of $Q_{1,0}$ and $Q_{0,1}$
at a particular point 
in the intersection 
$|z_{1}|>|z_{2}|>|z_{1}-z_{2}|>0$ as initial conditions for the 
system (\ref{system1})--(\ref{system4-1}). Then by the theory of 
differential equations of regular singular points, this initial value
problem has a unique solution of the form (\ref{ext}) in the 
region $|z_{2}|>|z_{1}-z_{2}|>0$
(without logarithm terms because the analytic 
extension of the solution in the region $|z_{1}|>|z_{2}|>0$
does not have logarithm terms). So $Q_{0, 1}$ and $Q_{1,0}$
can be 
analytically extended to analytic 
functions of the form (\ref{ext}) in the region $|z_2|>|z_1-z_2|>0$.

We have to do the same for $Q_{0,0}$ and $Q_{1,1}$.
By using the commutator formulas between $G(n+1/2)$, $L(n)$, 
$n\in \mathbb{Z}$, and the
intertwining operators, we have
\begin{eqnarray*}
 0&=&\langle w'_{(4)}, \mathcal{ Y}_{1}(w_{(1)}, z_{1})
\mathcal{ Y}_{2}(w_{(2)}, z_{2})Pw_{(3)}
\rangle\nn
&=&D_5(z_1,z_2)Q_{0,0}+D_6(z_1,z_2)Q_{1,1},\nn
0&=&\langle w'_{(4)}, \mathcal{ Y}_{1}(G(-1/2)w_{(1)}, z_{1})
\mathcal{ Y}_{2}(G(-1/2)w_{(2)}, z_{2})Pw_{(3)}
\rangle \nn
& =&D_7(z_1,z_2)Q_{0,0}+D_8(z_1,z_2)Q_{1,1} \nn
 0&=&\langle w'_{(4)}, \mathcal{ Y}_{1}(w_{(1)}, z_{1})
\mathcal{ Y}_{2}(Qw_{(2)}, z_{2})w_{(3)}
\rangle\nn
&=&D'_5(z_1,z_2)Q_{0,0}+D'_6(z_1,z_2)Q_{1,1},\nn
0&=&\langle w'_{(4)}, \mathcal{ Y}_{1}(G(-1/2)w_{(1)}, z_{1})
\mathcal{ Y}_{2}(G(-1/2)Qw_{(2)}, z_{2})w_{(3)}
\rangle \nn
& =&D'_7(z_1,z_2)Q_{0,0}+D'_8(z_1,z_2)Q_{1,1} 
\end{eqnarray*}
for some differential operators $D_5$, $D_6$, $D_7$, $D_8$, $D'_5$,
$D'_6$, $D'_7$ and $D'_8$. From the structures of $P$ and $Q$, we see
that these equations form a system of equations of regular singular
points with the only possible singular points $z_{1}, z_{2}=0, \infty$
and $z_{1}-z_{2}$. Using this system, we can prove the conclusions we
need for $Q_{0,0}$ and $Q_{1,1}$. We omit the details 
since they are completely the
same as those for $Q_{1, 0}$ and $Q_{0,1}$.
\epfv

The main goal of this section is the following more general result:

\begin{thm}\label{3-5}
For any quadruple $\{w_{(1)}, w_{(2)}, w_{(2)}, w'_{(4)}\}$ where 
$w_{(i)}\in W_{i}$, $i=1, 2, 3$, 
and $w'_{(4)}\in W'_{4}$ are homogeneous vectors,
the statement of 
Theorem \ref{top} is true. Moreover there is an integer $N$ depending only on
$\mathcal{Y}_{1}$ and $\mathcal{Y}_{1}$ such that 
\begin{equation}\label{wt-ineq}
\wt w_{(1)} +\wt w_{(2)} +s_{i}^{k, l}>N
\end{equation}
for all homogeneous elements $w_{(1)}\in W_{1}$, $w_{(2)}\in W_{2}$,
$i=1, \dots, m_{k, l}$ and $k, l=0, 1$. That is, in the terminology 
introduced in 
\cite{H1},
the products or the iterates of the intertwining  operators for
$L_{\mathfrak{n}\mathfrak{s}}(c_{p,q},0)$ have the convergence and extension property.
\end{thm}
\pf
We define the {\it weight} of a quadruple 
$\{w_{(1)}, w_{(2)}, w_{(2)}, w'_{(4)}\}$ to be 
$r=\wt w_{(1)}+\wt w_{(2)}+\wt w_{(3)}+\wt w'_{(4)}$.
We shall use induction on the weight $r$.
Theorem \ref{top} gives the result in the case of the smallest $r$. 

We have the following commutator formulas (which follow
{}from the Jacobi identity):
\bea 
\lefteqn{[Y(\omega,x_1),\mathcal{ Y}(w_{(1)},(x_2,\varphi))]}\nn
&& =(x_1^{-1}\delta(x_2/x_1)\pa_{x_2}+x_1^{-1}
\partial_{x_2}\delta(x_2/x_1)((\wt
w_{(1)})
+1/2\varphi \pa_{\varphi}))\mathcal{ Y}(w_{(1)},(x_2,\varphi))\nn
&&\quad +\cdots,
\label{formulas} \\
\lefteqn{[L(-n),\mathcal{ Y}(w_{(1)},(x_2,\varphi))]}\nn
&&= (x_2^{-n+1}\partial_{x_2}+(1-n)x_2^{-n}((\wt w_{(1)})
+1/2\varphi \pa_{\varphi}))\mathcal{ Y}(w_{(1)},(x_2,\varphi))\nn
&&\quad +\cdots,
\label{formulas-1} \\
\lefteqn{[Y(\tau,x_1), \mathcal{ Y}(w_{(1)},(x_2,\varphi))]} \nn
&&=(x_1^{-1}\delta(x_2/x_1)(\pa_{\varphi}-\varphi \pa_{x_2})-2
x_1^{-1}\partial_{x_2} \delta(x_2/x_1)(\wt w_{(1)})\varphi)
\mathcal{Y}(w_{(1)},(x_2,\varphi))\nn
&&\quad +\cdots,\label{formulas-2} \\
\lefteqn{[G(-n-1/2),
\mathcal{Y}(w_{(1)},(x_2,\varphi))]} \nn
&&= (x_2^{-n}(\partial_{\varphi}-\varphi \pa_{x_2})+2n
x_2^{-n-1}(\wt w_{(1)})
\varphi)\mathcal{ Y}(w_{(1)},(x_2,\varphi))+\cdots
\label{formulas-3}
\eea
where we write $\cdots$ for terms  associated to elements whose weights are less
than the weight of $w_{(1)}$.

To prove the theorem it is enough
to show that if the statement is true for the quadruple
$\{w_{(1)}, w_{(2)},w_{(3)},w'_{(4)}\}$, 
then it  is also true for  the quadruples
$$\{G(-n-1/2)w_{(1)}, w_{(2)},w_{(3)},w'_{(4)}\}, 
\{L(-n)w_{(1)}, w_{(2)},w_{(3)},w'_{(4)}\},$$
$$\{w_{(1)}, G(-n-1/2)w_{(2)},w_{(3)},w'_{(4)}\}, 
\{w_{(1)}, L(-n)w_{(2)},w_{(3)},w'_{(4)}\}, $$
$$\{w_{(1)}, w_{(2)},G(-n-1/2)w_{(3)},w'_{(4)}\}, \{w_{(1)}, w_{(2)},
L(-n)w_{(3)},w'_{(4)}
\},$$ 
$$\{w_{(1)}, w_{(2)},w_{(3)},G(-n-1/2)w'_{(4)}\}, 
\{w_{(1)}, w_{(2)},w_{(3)},L(-n)w'_{(4)}\}$$
for every $n \in \mathbb{N}$. 
We shall prove the statement only for  the quadruples
$$\{w_{(1)}, w_{(2)},G(-n-1/2)w_{(3)},w'_{(4)}\}, \{w_{(1)}, w_{(2)},
L(-n)w_{(3)},w'_{(4)}
\};$$ 
for the others, the proofs are similar.

{}From (\ref{formulas})--(\ref{formulas-3}) we get
\bea
\lefteqn{\langle w'_{(4)}, \mathcal{ Y}_{1}(w_{(1)}, (x_1,\varphi_1))
\mathcal{Y}_{2}(w_{(2)}, (x_2,\varphi_2))L(-n)w_{(3)} \rangle} \nn
&&=-(x_1^{-n+1}\partial_{x_1}+(1-n)x_1^{-n}((\wt w_{(1)})
+1/2\varphi_1 \pa_{\varphi_1})+ x_2^{-n+1}\partial_{x_2} \nn
&&\quad  +(1-n)x_2^{-n}((\wt w_{(2)})+1/2\varphi_2 \pa_{\varphi_2})
\cdot \nn
&&\quad\quad\cdot
\langle w'_{(4)}, \mathcal{ Y}_{1}(w_{(1)},
(x_1,\varphi_1))\mathcal{Y}_{2}(w_{(2)}, (x_2,\varphi_2))w_{(3)} 
\rangle \nn
&&\quad  +\langle L(n) w'_{(4)}, \mathcal{Y}_{1}(w_{(1)},
(x_1,\varphi_1))\mathcal{Y}_{2}(w_{(2)}, (x_2,\varphi_2))w_{(3)} \rangle
+\cdots
\eea
and
\bea
\lefteqn{\langle w'_{(4)}, \mathcal{ Y}_{1}(w_{(1)}, (x_1,\varphi_1))
\mathcal{Y}_{2}(w_{(2)}, (x_2,\varphi_2))G(-n-1/2)w_{(3)} \rangle} \nn
&&= -(x_1^{-n}(\partial_{\varphi_1}-\varphi_1 \pa_{x_1})+2n
x_1^{-n-1}(\wt w_{(1)}))\varphi_1+x_2^{-n}(\partial_{\varphi_2} \nn
&&\quad - \varphi_2 \pa_{x_2})+2n x_2^{-n-1}(\wt w_{(2)})\varphi_2)
\cdot \nn 
&&\quad\quad \cdot\langle w'_{(4)}, \mathcal{ Y}_{1}(w_{(1)},
(x_1,\varphi_1))\mathcal{ Y}_{2}(w_{(2)}, (x_2,\varphi_2))w_{(3)} \rangle \nn
&& \quad +\langle G(n+1/2) w'_{(4)}, \mathcal{ Y}_{1}(w_{(1)},
(x_1,\varphi_1))\mathcal{ Y}_{2}(w_{(2)}, (x_2,\varphi_2))w_{(3)} 
\rangle \nn
&& \quad +\cdots,
\eea
where $\cdots$ 
represents terms 
(certain derivatives of matrix coefficients) involving quadruples 
whose weights are less
then $r$.
Since the weights of quadruples $\{w_{(1)},w_{(2)},w_{(3)}, 
L(n) w'_{(4)} \}$ and
$\{w_1,w_2,w_3, G(n+1/2) w'_{4} \}$ are $r-n$ and $r-n-1/2$
respectively we can apply induction. 

We choose $N$ to be an integer such that when $w_{(1)}$ and $w_{(2)}$ are 
the lowest weight vectors, (\ref{wt-ineq}) holds. Then using the formulas 
above and induction, it is easy to see that (\ref{wt-ineq}) holds 
for all $w_{(1)}\in W_{1}$ and $w_{(2)}\in W_{2}$ (cf. 
\cite{H2}).  
\epfv

Finally, we generalize Theorem \ref{3-5} to vertex operator algebras
in the class $\mathcal{C}_{p_{1}, q_{1};
\dots; p_{n}, q_{n}}$ for $(p_{i}, q_{i})$, 
$i=1, \dots, n$, $n$
pairs of integers larger than $1$ such that $p_{i}-q_{i}\in
2\mathbb{Z}$ and $(p_{i}-q_{i})/2$ and $q_{i}$ are relatively prime to
each other. Since the proof of the following 
 result is completely the same as the
corresponding results in \cite{H2} and \cite{HL5.5}, 
we omit the details
here.

\begin{thm}\label{3-6}
Let $V$ be a vertex operator algebra in the class $\mathcal{C}_{p_{1},
q_{1}; \dots; p_{n}, q_{n}}$. Then the products and iterates of
intertwining operators for $V$ have the convergence and extension
property. \epf
\end{thm}

\renewcommand{\theequation}{\thesection.\arabic{equation}}
\renewcommand{\thethm}{\thesection.\arabic{thm}}
\setcounter{equation}{0}
\setcounter{thm}{0}

\section{Intertwining operator superalgebras and 
vertex tensor categories for $N=1$ minimal models}

In this section, $V$ is a vertex operator algebra in the class 
$\mathcal{C}^{(1)}_{p_{1}, q_{1}; \dots;  
p_{n}, q_{n}}$ for $(p_{i}, q_{i})$, 
$i=1, \dots, n$, $n$
pairs of integers larger than $1$ such that $p_{i}-q_{i}\in
2\mathbb{Z}$ and $(p_{i}-q_{i})/2$ and $q_{i}$ are relatively prime to
each other.  
By Proposition
\ref{2-7} and Theorem \ref{3-6}, and Theorems 3.1 and 3.2
in \cite{H2}, which are in turn proved in 
 \cite{H2} using results
in \cite{HL1}--\cite{HL5} and \cite{H1}, we obtain the
following:

\begin{thm}[associativity for intertwining operators]

\begin{enumerate}

\item For any $V$-modules 
$W_{0}$, $W_{1}$, $W_{2}$, $W_{3}$  and $W_{4}$, 
any 
intertwining operators $\mathcal{Y}_{1}$ and 
$\mathcal{Y}_{2}$ of 
 types ${W_{0}}\choose {W_{1}W_{4}}$ and ${W_{4}}
\choose {W_{2}W_{3}}$,
respectively,
and any choice of $\log z_{1}$ and $\log z_{2}$,
$$\langle w'_{(0)}, \mathcal{Y}_{1}(w_{(1)},
x_{1})\mathcal{Y}_{2}(w_{(2)}, x_{2})w_{(3)}\rangle\lbar_{x^{n}_{1}
=e^{n\log z_{1}},\;
x^{n}_{2}=e^{n\log z_{2}}, \; n\in \mathbb{C}}$$
is absolutely convergent when $|z_{1}|>|z_{2}|>0$ for
$w'_{(0)}\in W'_{0}$, $w_{(1)}\in W_{1}$,
$w_{(2)}\in W_{2}$ and $w_{(3)}\in W_{3}$.
For any modules 
$W_{0}$, $W_{1}$, $W_{2}$, $W_{3}$,  and $W_{5}$ and 
any intertwining operators $\mathcal{Y}_{3}$
and $\mathcal{Y}_{4}$ of types ${W_{5}}\choose {W_{1}W_{2}}$ and 
${W_{0}}\choose {W_{5}W_{3}}$, respectively, and any choice of 
$\log z_{2}$ and $\log (z_{1}-z_{2})$,
$$\langle w'_{(0)}, \mathcal{Y}_{4}(\mathcal{Y}_{3}(w_{(1)},
x_{0})w_{(2)}, x_{2})w_{(3)}\rangle\lbar_{x^{n}_{0}=e^{n \log
(z_{1}-z_{2})},\;
x^{n}_{2}=e^{n\log z_{2}}, \; n\in \mathbb{C}}$$
is absolutely convergent when $|z_{2}|>|z_{1}-z_{2}|>0$
for $w'_{(0)}\in W'_{0}$, $w_{(1)}\in W_{1}$,
$w_{(2)}\in W_{2}$ and $w_{(3)}\in W_{3}$.

\item For any $V$-modules 
$W_{0}$, $W_{1}$, $W_{2}$, $W_{3}$  and $W_{4}$,
any 
intertwining operators $\mathcal{Y}_{1}$ and $\mathcal{Y}_{2}$ of 
types ${W_{0}}\choose {W_{1}W_{4}}$ and ${W_{4}}\choose {W_{2}W_{3}}$,
respectively, 
there exist a module $W_{5}$ and intertwining operators 
$\mathcal{Y}_{3}$
and $\mathcal{Y}_{4}$ of types ${W_{5}}\choose {W_{1}W_{2}}$ and 
${W_{0}}\choose {W_{5}W_{3}}$, respectively, such that for
any $z_{1}, z_{2}\in \mathbb{C}$ satisfying 
$|z_{1}|>|z_{2}|>|z_{1}-z_{2}|>0$ and for any 
$w'_{(0)}\in W'_{0}$, $w_{(1)}\in W_{1}$,
$w_{(2)}\in W_{2}$ and $w_{(3)}\in W_{3}$,
\begin{eqnarray}\label{1-1}
\lefteqn{\langle w'_{(0)}, \mathcal{Y}_{1}(w_{(1)},
x_{1})\mathcal{Y}_{2}(w_{(2)}, x_{2})w_{(3)}
\rangle|_{x_{1}^{n}=e^{n\log z_{1}}, x_{2}^{n}=e^{n\log z_{2}},
n\in \mathbb{R}}}
\nn
&&=\langle w'_{(0)}, \mathcal{Y}_{4}(\mathcal{Y}_{3}(w_{(1)},
x_{0})w_{(2)}, x_{2})w_{(3)}\rangle
\rangle|_{x_{0}^{n}=e^{n\log (z_{1}-z_{2})}, 
x_{2}^{n}=e^{n\log z_{2}},
n\in \mathbb{R}},\nn
&&
\end{eqnarray}
where $\log z_{1}=|z_{1}|+i\arg z_{1}$, $\log z_{2}
=|z_{2}|+i\arg z_{2}$ and $\log (z_{1}-z_{2})=
|z_{1}-z_{2}|+i\arg (z_{1}-z_{2})$ are 
the values of the logarithms of $z_{1}$, $z_{2}$ and $z_{1}-z_{2}$
such that $0\le \arg z_{1}, \arg z_{2}, 
\arg (z_{1}-z_{2})\le 2\pi$.

\item For any modules 
$W_{0}$, $W_{1}$, $W_{2}$, $W_{3}$,  and $W_{5}$, 
any intertwining operators $\mathcal{Y}_{3}$
and $\mathcal{Y}_{4}$ of types ${W_{5}}\choose {W_{1}W_{2}}$ and 
${W_{0}}\choose {W_{5}W_{3}}$, respectively, 
there exist a module $W_{4}$ and 
intertwining operators $\mathcal{Y}_{1}$ and $\mathcal{Y}_{2}$ of 
types ${W_{0}}\choose {W_{1}W_{4}}$ and ${W_{4}}\choose {W_{2}W_{3}}$,
respectively, such that for any $z_{1}, z_{2}\in \mathbb{C}$ 
satisfying 
$|z_{1}|>|z_{2}|>|z_{1}-z_{2}|>0$ and for any
$w'_{(0)}\in W'_{0}$, $w_{(1)}\in W_{1}$,
$w_{(2)}\in W_{2}$ and $w_{(3)}\in W_{3}$, the equality
(\ref{1-1})  holds.\epf

\end{enumerate}
\end{thm}

\begin{thm}[commutativity for intertwining operators]
For any $V$-modules 
$W_{0}$, $W_{1}$, $W_{2}$, $W_{3}$  and $W_{4}$ and 
any 
intertwining operators $\mathcal{Y}_{1}$ and $\mathcal{Y}_{2}$ of 
types ${W_{0}}\choose {W_{1}W_{4}}$ and ${W_{4}}\choose {W_{2}W_{3}}$,
respectively,
there exist a module $W_{5}$ and intertwining operators $\mathcal{Y}_{3}$
and $\mathcal{Y}_{4}$ of types ${W_{0}}\choose {W_{2}W_{5}}$ and 
${W_{5}}\choose {W_{1}W_{3}}$, respectively, such that for any 
homogeneous $w'_{(0)}\in W'_{0}$, $w_{(1)}\in W_{1}$,
$w_{(2)}\in W_{2}$ and $w_{(3)}\in W_{3}$,
the multivalued analytic function 
$$\langle w'_{(0)}, \mathcal{Y}_{1}(w_{(1)},
x_{1})\mathcal{Y}_{2}(w_{(2)}, x_{2})w_{(3)}\rangle\lbar_{x_{1}
=z_{1},\;
x_{2}=z_{2}}$$
of $z_{1}$ and $z_{2}$ in the region $|z_{1}|>|z_{2}|>0$ and 
the multivalued analytic function 
$$(-1)^{|w_{(1)}||w_{(2)}|}\langle w'_{(0)}, \mathcal{Y}_{3}(w_{(2)},
x_{2})\mathcal{Y}_{4}(w_{(1)}, x_{1})w_{(3)}\rangle\lbar_{x_{1}
=z_{1},\;
x_{2}=z_{2}}$$
of $z_{1}$ and $z_{2}$ in the region $|z_{2}|>|z_{1}|>0$ are analytic
extensions of each other.\epf
\end{thm}

In \cite{H2.5}, the notion of intertwining operator algebra
was introduced (see also \cite{H4} and \cite{H6}). 
This notion has the following generalization:

\begin{defn}
{\it An $\mathbb{R}\times \mathbb{Z}_{2}$-graded 
vector space $W=\coprod_{n\in \mathbb{R}} 
W_{(n)}$, together with a finite set $\mathcal{A}$ with a
distinguished element $e$, a subspace 
$\mathcal{V}_{a_{1}a_{2}}^{a_{3}}$ of the space of 
linear maps {}from $W^{a_{1}}\otimes W^{a_{2}}\to 
W^{a_{3}}\{x\}$ for $a_{1}, a_{2}, a_{3}\in \mathcal{A}$, 
two distinguished vectors, the {\it vacuum} $\mathbf{1}$ and 
the {\it Virasoro element} $\omega$ of $W^{e}$, 
is called an {\it intertwining 
operator superalgebra} if $W$ 
together with the other data satisfies 
all the axioms for intertwining operator algebras except that 
commutativity (or skew-symmetry) for intertwining operators
is replaced by the 
corresponding commutativity  (or skew-symmetry)
for intertwining operators 
for vertex operator superalgebras. Note that 
there is a unique intertwining operator $Y$ of type ${W^{e}\choose
W^{e}W^{e}}$ such that $(W^{e}, Y, \mathbf{1}, \omega)$ is a vertex
operator subalgebra of $W$. If there is an element 
$\tau\in W_{e}$ such that $(V, Y, \mathbf{1}, \tau)$  is 
an $N=1$ superconformal vertex operator algebra, then $W$ together
with all the other data is
called an {\it $N=1$ superconformal intertwining 
operator superalgebra}.}
\end{defn}

The intertwining operator superalgebra just defined is denoted 
$$(W, \mathcal{A}, \{\mathcal{V}_{a_{1}a_{2}}^{a_{3}}\}, 
\mathbf{1}, \tau)$$
or simply $W$. 
Theorem 3.5
in \cite{H2.5} can be modified easily to incorporate the 
signs and to apply to the cases discussed in the present paper. 
Applying this modification
to $V$, we obtain the following:

\begin{thm}
Assume in addition that $V$ is rational. 
Let $\mathcal{A}=\{a_{i}\}_{i=1}^{m}$ 
be  the set of all equivalence 
classes of
irreducible $V$-modules.  Let $W^{a_{1}}, \dots, W^{a_{m}}$ be
representatives of $a_{1}, \dots, a_{m}$, respectively.  Let
$W=\coprod_{i=1}^{m}W^{a_{i}}$, and let $\mathcal{V}_{a_{1}a_{2}}^{a_{3}}$, 
for $a_{1}, a_{2}, 
a_{3}\in \mathcal{A}$, be the space
of intertwining operators of type ${W^{a_{3}}\choose
W^{a_{1}}W^{a_{2}}}$. Then $(W, \mathcal{A}, 
\{\mathcal{V}_{a_{1}a_{2}}^{a_{3}}\}, \mathbf{1}, \tau)$ 
(where $\mathbf{1}$,
$\tau$ are the vacuum and Neuve-Schwarz
element of $V$, respectively) is an $N=1$ superconformal 
intertwining operator superalgebra.\epf
\end{thm}

In particular, we have:

\begin{thm}
For any integers $p, q$  larger than $1$ such that 
$p-q\in 2\mathbb{Z}$ and $(p-q)/2$ and $q$ 
 relatively prime to each other, the 
direct sum 
$$\coprod_{m, n\in \mathbb{Z}, m-n\in 2\mathbb{Z},
0<m<p, 0<n<q}
L_{\mathfrak{n}\mathfrak{s}}(c_{p, q}, h_{p, q}^{m, n})$$
together with the finite set 
$$\{m, n\in \mathbb{Z}, m-n\in 2\mathbb{Z},
0<m<p, 0<n<q\},$$
the spaces of intertwining operators of type
$${L_{\mathfrak{n}\mathfrak{s}}(c_{p, q}, h_{p, q}^{m_{3}, n_{3}})\choose
L_{\mathfrak{n}\mathfrak{s}}(c_{p, q}, h_{p, q}^{m_{1}, n_{1}})
L_{\mathfrak{n}\mathfrak{s}}(c_{p, q}, h_{p, q}^{m_{2}, n_{2}})}$$
for $m_{i}, n_{i}\in \mathbb{Z}, m_{i}-n_{i}\in 2\mathbb{Z},
0<m_{i}<p, 0<n_{i}<q$, $i=1, 2$, and the vacuum and 
the Neveu-Schwarz elements 
of $L_{\mathfrak{n}\mathfrak{s}}(c_{p, q}, 0)$
is an $N=1$ superconformal  intertwining operator superalgebra.\epf
\end{thm}

Now we discuss the vertex tensor category structures. 
Recall the sphere partial operad $K=\{K(j)\}_{j\in \mathbb{N}}$, the
vertex partial operads $\tilde{K}^{c}=\{\tilde{K}^{c}(j)\}_{j\in \mathbb{N}}$ 
of central
charge $c\in \mathbb{C}$ constructed in \cite{H3} and 
the definition of vertex tensor category in \cite{HL4} and 
\cite{HL6}. For any $c\in \mathbb{C}$ and $j\in \mathbb{N}$, $\tilde{K}^{c}(j)$
is a trivial holomorphic line bundle over $K(j)$ and we have a 
canonical holomorphic section $\psi_{j}$. Given a vertex tensor category,
we have, among other things, 
a tensor product bifunctor $\boxtimes_{\tilde{Q}}$ for each
$\tilde{Q}\in \tilde{K}^{c}(2)$. In particular,  
$\psi_{2}(P(z))\in \tilde{K}^{c}(2)$ and thus there is a 
tensor product bifunctor $\boxtimes_{\psi_{2}(P(z))}$.

Note that $\boxtimes_{P(z)}$ constructed
in \cite{HL5} can be generalized without any difficulty 
to categories of modules for vertex operator superalgebras. 
By Proposition \ref{2-7} and Theorem \ref{3-6}, and
Theorem 3.2 and Corollary 3.3
in \cite{H2}, we obtain:

\begin{thm}\label{vtc}
Let $c$ be the central charge of $V$.
Then the category of $V$-modules has a natural structure of  vertex
tensor category of central charge $c$ such that for each $z\in 
\mathbb{C}^{\times}$, the tensor product bifunctor $\boxtimes_{\psi_{2}(P(z))}$
associated with $\psi_{2}(P(z))\in
\tilde{K}^{c}(2)$ is equal to generalization to the category $V$-modules 
of 
$\boxtimes_{P(z)}$ constructed
in \cite{HL5}.\epf
\end{thm}

Combining Theorem \ref{vtc}  with Theorem 4.4 in \cite{HL4} (see
\cite{HL6} for the proof), we obtain:

\begin{cor}
The category of $V$-modules has a natural 
structure of  braided tensor category such that the tensor product
bifunctor is $\boxtimes_{P(1)}$. In particular, 
the category of $L_{\mathfrak{n}\mathfrak{s}}(c_{p_{1}, q_{1}}, 0)\otimes \cdots \otimes 
L_{\mathfrak{n}\mathfrak{s}}(c_{p_{m}, q_{m}}, 0)$-modules has a natural 
structure of   braided tensor category.\epf
\end{cor}

In particular, the special case $V=L_{\mathfrak{n}\mathfrak{s}}(c_{p, q}, 0)$ gives:

\begin{thm}
For any integers $p, q$  larger than $1$ such that 
$p-q\in 2\mathbb{Z}$ and $(p-q)/2$ and $q$ 
 relatively prime to each other,
the category of modules for the $N=1$ Neveu-Schwarz Lie 
superalgebra equivalent to
$L_{\mathfrak{n}\mathfrak{s}}(c_{p, q}, h_{p, q}^{m, n})$, $m, n\in \mathbb{Z}$, $m-n\in 2\mathbb{Z}$,
$0<m<p, 0<n<q$, has a natural 
structure of  braided tensor category such that the tensor product
bifunctor is $\boxtimes_{P(1)}$.\epf
\end{thm}

\renewcommand{\theequation}{\thesection.\arabic{equation}}
\renewcommand{\thethm}{\thesection.\arabic{thm}}
\setcounter{equation}{0}
\setcounter{thm}{0}

\section{Appendix: An example}

In this appendix we give a concrete example.
Using this example, we present a way to calculate the
fusion rules for minimal models and also show how to
derive explicit differential equations for
matrix coefficients.

By Theorem \ref{ad}, 
for integers $p, q$  larger than $1$ such that 
$p-q\in 2\mathbb{Z}$ and $(p-q)/2$ and $q$ are 
 relatively prime to each other, all irreducible
representations of vertex operator superalgebra
$L_{\mathfrak{n}\mathfrak{s}}(c_{p,q},0)$ are given by the
set $L_{\mathfrak{n}\mathfrak{s}}(c_{p,q}, h_{p, q}^{m,n})$ 
where $m, n \in \mathbb{N}$, 
$0<m<p,0<n<q$ and
$m-n \in 2\mathbb{Z}$.
We consider the special case $c=-11/14$, that is, $p=7$ and $q=3$.
Then vertex operator superalgebra
$L_{\mathfrak{n}\mathfrak{s}}(-11/14,0)$ has three inequivalent irreducible modules
$L_{\mathfrak{n}\mathfrak{s}}(-11/14, 2/7)$, $L_{\mathfrak{n}\mathfrak{s}}(-11/14, -1/14)$ and itself. 
The following fact follows {}from the representation 
theory (see \cite{As}) of $N=1$ Neveu-Schwarz Lie superalgebra:

\begin{prop} \label{ex1}
\begin{enumerate}

\item The irreducible module $L_{\mathfrak{n}\mathfrak{s}}(-11/14,0)$ 
is equivalent to $V_{\mathfrak{ns}}(-11/14,0)/J^1$, where
$J^1$ is a submodule generated by a (singular)
vector of weight $6$. 

\item The irreducible module
$L_{\mathfrak{n}\mathfrak{s}}(-11/14,2/7)$ is equivalent to 
$$M_{\mathfrak{n}\mathfrak{s}}(-11/14,2/7)/(J^2+J^3)$$
where $J^2$  and $J^3$ are submodules generated by
(singular) vectors of weights $2$  and $5/2$, respectively.

\item The irreducible module
$L_{\mathfrak{n}\mathfrak{s}}(-11/14,-1/14)$ is equivalent to 
$$M_{\mathfrak{n}\mathfrak{s}}(-11,14,-1/14)/(J^4+J^5)$$
where  $J^4$  and $J^5$ are submodules generated by
(singular) vectors of weights $4$ and $3/2$, respectively.\epf

\end{enumerate}
\end{prop}

The proof of this result is easy and is omitted.

Now we assume that the reader is familiar with 
the theory of Zhu's algebra associated to a vertex operator
algebra (see \cite{Z}). 
It has been shown in \cite{KW} that  $A(V_{\mathfrak{ns}}(c,0))$
(Zhu's algebra associated to 
$V_{\mathfrak{ns}}(c,0)$) is
isomorphic to $\mathbb{C}[x]$ where $x=[\omega]$. 
Using this information
and Conclusion 1 in Proposition \ref{ex1}, we see that
$A(L_{\mathfrak{n}\mathfrak{s}}(-11/14,0))$ is isomorphic to 
$\mathbb{C}[x]/\langle x(x-2/7)(x+1/14)\rangle$.

Next we use Frenkel-Zhu's formula in \cite{FZ} 
to calculate the fusion rules. 
First we have the following  result which can be obtained by using
the results in  \cite{KW}:

\begin{prop}\label{5-2}
Let $A(M_{\mathfrak{n}\mathfrak{s}}(c,h))$ be the $A(V_{\mathfrak{ns}}(c,0))$-bimodule associated to the 
$V_{\mathfrak{ns}}(c, 0)$-module $M_{\mathfrak{n}\mathfrak{s}}(c, h)$ with
the left  and right action given by
\bea
a*m&=&\res_x Y(a,x)\frac{(1+x)^{\deg(a)}}{x}m,\\
m*a&=&\res_x Y(a,x)\frac{(1+x)^{\deg(a)-1}}{x}m,
\eea
respectively, 
for $a \in V_0$, $m\in A(M_{\mathfrak{n}\mathfrak{s}}(c, h))$
and with trivial action for $a \in V_1$, 
$m\in A(M_{\mathfrak{n}\mathfrak{s}}(c, h))$.
Then we have:
\begin{enumerate}
\item The $A(V_{\mathfrak{ns}}(c,0))$-bimodule
$A(M_{\mathfrak{n}\mathfrak{s}}(c,h))$ is equivalent to 
$$\mathbb{C}[x,y] 
\oplus \mathbb{C}[x,y]v,$$
where $x=[L(0)+2L(-1)+L(-2)]$, $y=[L(-1)+L(-2)]$ and
$v=[G(-1/2)]$.

\item
We have the following action of $A(V_{\mathfrak{n}\mathfrak{s}}(c,0))\cong \mathbb{C}[y]$
on $A(M_{\mathfrak{n}\mathfrak{s}}(c,h))\cong \mathbb{C}[x,y] 
\oplus \mathbb{C}[x,y]v$:
$$y*(x^k y^l)=x^{k+1} y^l, \ \ (x^k y^l)*y=x^k y^{l+1},$$
$$y*(x^k y^lv)=x^{k+1} y^lv, \ \ (x^k y^lv)*y=x^k y^{l+1}v.\epfe$$

\end{enumerate}
\end{prop}

The proof of this result is straightforward and is omitted.

{}From this result, we see that $A(M_{\mathfrak{n}\mathfrak{s}}(c,h))$ has ${\bf Z}_2$-grading:
$$A(M_{\mathfrak{n}\mathfrak{s}}(c,h))=A(M_{\mathfrak{n}\mathfrak{s}}(c,h))^{+} \oplus A(M_{\mathfrak{n}\mathfrak{s}}(c,h))^{-}$$ 
where 
$$A(M_{\mathfrak{n}\mathfrak{s}}(c,h))^{+}=\mathbb{C}[x,y]$$ and 
$$A(M_{\mathfrak{n}\mathfrak{s}}(c,h))^{-}=\mathbb{C}[x,y]v.$$

Using this result, we identify
$A(L_{\mathfrak{n}\mathfrak{s}}(-11/14,2/7))$ as follows: It is not hard to see that
$J^2$ is generated by the vector 
$$(8L(-2)-21L(-1)^2+21G(-3/2)G(-1/2))w$$
and $J^3$ is generated by
$$\left(L(-2)G(-1/2)-\frac{25}{14}
G(-5/2)-\frac{7}{6}L(-1)^2G(-1/2)+\frac{3}{2}L(-1)G(-3/2)\right)w,$$
where $w$ is a lowest weight vector in $M_{\mathfrak{n}\mathfrak{s}}(-11/14,2/7)$. 

Since for any $A(L_{\mathfrak{n}\mathfrak{s}}(-11/14,0))$-module $M$ and 
any submodule $I$ of $M$, 
$A(M/I)$ is equivalent to $A(M)/A(I)$, we obtain 
by using Proposition \ref{5-2} and  calculations that 
$A(L_{\mathfrak{n}\mathfrak{s}}(-11/14,2/7))$ is equivalent to
$\mathbb{C}[x,y]/I_1 \oplus \mathbb{C}[x,y]v/I_2v$ where
\begin{eqnarray*}
I_1&=&\left\langle -21(x-y)^2+4(x+y)+\frac{4}{7}, 
(x-y)(49(x-y)^2-84(x+y)+20)\right\rangle, \\
I_2&=&\left\langle -\frac{7}{6}(x-y)^2+\frac{1}{2}(x+y)+\frac{1}{24}, 
-\frac{75}{28}+25(x+y)-21(x-y)^2\right\rangle.
\end{eqnarray*}
Note that $\dim_{\mathbb{C}} A(L_{\mathfrak{n}\mathfrak{s}}(-11/14,2/7))=4$.

We now apply Frenkel-Zhu's formula for the fusion rules, that is,
the
dimensions of the spaces of intertwining operators. 
It is proved in
\cite{FZ} (see also \cite{L})
that for  suitable vertex operator algebra $V$ and any $V$-modules 
$W_{1}$, $W_{2}$ and $W_{3}$,  the fusion rule
$$\mathcal{N}_{W_{1} W_{2}}^{W_{3}}=\dim \hom_{A(V)} 
(A(W_{1})\otimes _{A(V)}W_{2}(0), W_{3}(0)),$$
where $W_{i}(0)$, $i=1, 2, 3$,
are the top level of $W_{i}$, respectively,
equipped with $A(V)$-module structures.

We have shown in the proof of Proposition \ref{fusion}
 that any intertwining operator 
$\mathcal{ Y}$ among irreducible modules
is uniquely determined by $\mathcal{ Y}(w_{1},z)$ and
$\mathcal{ Y}(G(-1/2)w_{1},z)$ where 
$w_{(1)}$ is a lowest weight vector.
Now we have actually reproved this fact
using Frenkel-Zhu's theory. In fact {}from the $\mathbb{Z}_2$-grading
of $A(W)=A(W)^{0} \oplus A(W)^{1}$,  it follows that 
\begin{eqnarray*}
\lefteqn{\hom_{A(V)} (A(W_{1})\otimes _{A(V)}W_{2}(0), W_{3}(0)) 
\cong}\nn
&&\hom_{A(V)} (A(W_{1})^{0} \otimes _{A(V)}W_{2}(0), W_{3}(0)) \nn
&&\quad\quad \oplus
\hom_{A(V)} (A(W_{1})^{-}\otimes _{A(V)}W_{2}(0), W_{3}(0)).
\end{eqnarray*}
Thus one can define in an 
obvious way $(\mathcal{N}_{W_{1}W_{2}}^{W_{3}})^{0}$ and 
$(\mathcal{N}_{W_{1}W_{2}}^{W_{3}})^{1}$. 
We have proved in Proposition \ref{fusion}
 that $\mathcal{N}_{W_{1}W_{2}}^{W_{3}}\leq 2$.
A natural question is: Is it true that 
$(\mathcal{N}_{W_{1}W_{2}}^{W_{3}})^{0}
\leq 1$ and $(\mathcal{N}_{W_{1}W_{2}}^{W_{3}})^{1}\leq 1$?

Let us go back to our example. For convenience, we shall 
denote the fusion rules among irreducible 
modules
$L_{\mathfrak{n}\mathfrak{s}}(-11/14, h_{1})$, 
$L_{\mathfrak{n}\mathfrak{s}}(-11/14, h_{2})$ 
and $L_{\mathfrak{n}\mathfrak{s}}(-11/14, h_{3})$
for $L_{\mathfrak{n}\mathfrak{s}}(-11/14,0)$ by 
$\mathcal{ N}_{h_{1} h_{2}}^{h_{3}}$.
By using the formula above and this notation, we obtain
\begin{eqnarray*}
\mathcal{ N}_{2/7,2/7}^{0}&=&1, \\
\mathcal{ N}_{2/7,2/7}^{-1/14}&=&1, \\
\mathcal{ N}_{2/7,2/7}^{2/7}&=&0, \\
\mathcal{N}_{2/7,-1/14}^{0}&=&0, \\
\mathcal{ N}_{2/7,-1/14}^{2/7}&=&1,\\
\mathcal{ N}_{2/7, -1/14}^{2/7}&=&1,\\
\mathcal{ N}_{2/7,-1/14}^{-1/14}&=&1,\\
\mathcal{ N}_{2/7, -1/14}^{-1/14}&=&1,
\end{eqnarray*}
in addition to the obvious fusion rules.
In the proceeding section, we have 
proved that there is a braided tensor
category structure on the category of modules for
$L_{\mathfrak{n}\mathfrak{s}}(c_{p,q},0)$. Using the tensor product notation
$\boxtimes=\boxtimes_{P(1)}$ and using $\cong$ as an abbreviation of
the phrase ``is isomorphic 
to,'' we have the following 
decompositions of the tensor product modules:
\begin{eqnarray*}
\lefteqn{L_{\mathfrak{n}\mathfrak{s}}(-11/14,2/7) \boxtimes L_{\mathfrak{n}\mathfrak{s}}(-11/14,2/7)}\nn
&&\cong
L_{\mathfrak{n}\mathfrak{s}}(-11/14,-1/14)\oplus L_{\mathfrak{n}\mathfrak{s}}(-11/14,0),\\
\lefteqn{L_{\mathfrak{n}\mathfrak{s}}(-11/14,2/7) \boxtimes L_{\mathfrak{n}\mathfrak{s}}(-11/14,-1/14)}\nn
&&\cong 
L_{\mathfrak{n}\mathfrak{s}}(-11/14,2/7)\oplus L_{\mathfrak{n}\mathfrak{s}}(-11/14,-1/14),\\
\lefteqn{L_{\mathfrak{n}\mathfrak{s}}(-11/14,-1/14) \boxtimes  L_{\mathfrak{n}\mathfrak{s}}(-11/14,-1/14)}\nn
&&\cong
L_{\mathfrak{n}\mathfrak{s}}(-11/14,0) \oplus L_{\mathfrak{n}\mathfrak{s}}(-11/14,-1/14) \oplus L_{\mathfrak{n}\mathfrak{s}}(-11/14,2/7).
\end{eqnarray*}

Finally we discuss the differential equations.  Let 
$$\mathcal{Y}_1,
\mathcal{Y}_2 \in \mathcal{V}^{L_{\mathfrak{n}\mathfrak{s}}(-11/14,2/7)}
_{L_{\mathfrak{n}\mathfrak{s}}(-11/14,-1/14)\;\;
L_{\mathfrak{n}\mathfrak{s}}(-11/14,2/7)}.$$  
Our  goal is to 
derive differential equations {}from which we
can find the sums of the series
$$\langle w'_{(4)},
\mathcal{ Y}_{1}(w_{(1)}, z_{1}) \mathcal{ Y}_{2}(w_{(2)},
z_{2})w_{(3)} \rangle, 
$$
where $w_{(1)}, w_{(2)} \in
L_{\mathfrak{n}\mathfrak{s}}(-11/14,-1/14)$, $w_{(3)}\in L_{\mathfrak{n}\mathfrak{s}}(-11/14,2/7)$ and $w'_{(4)}\in
L'(-11/14,2/7)$ are the lowest weight vectors in these modules.  
We will see that even with  singular vectors of
weight $2$, we have highly non-trivial differential equations.

Since
\begin{eqnarray*}
\lefteqn{\langle w'_{(4)}, \mathcal{ Y}_{1}(w_{(1)}, z_{1})
\mathcal{ Y}_{2}(w_{(2)}, z_{2})(8L(-2)-21L(-1)^2} \nn
&& \quad\quad\quad\quad\quad\quad\quad\quad\quad \quad\quad\quad
+21G(-3/2)G(-1/2))w_{(3)} \rangle=0, \nn
\lefteqn{\langle w'_{(4)}, \mathcal{ Y}_{1}(G(-1/2)w_{(1)}, z_{1})
\mathcal{ Y}_{2}(G(-1/2)w_{(2)}, z_{2})(8L(-2)-21L(-1)^2} \nn
&& \quad\quad\quad\quad\quad\quad\quad\quad\quad \quad\quad\quad
+21G(-3/2)G(-1/2))w_{(3)} \rangle=0,
\end{eqnarray*}
 we obtain 
\begin{eqnarray*}
\lefteqn{\left(-8\left(\frac{1}{14z_1^{2}}
+\frac{1}{14z_2^{2}}\right)-
21(\partial_{z_1}+\partial_{z_2})^2-29
\left(\frac{\partial_{z_1}}{z_1}
+\frac{\partial_{z_2}}{z_2}\right)\right)Q_{0,0}(z_1,z_2)} \nn
&&\quad\quad\quad\quad\quad\quad \quad\quad\quad\quad \quad\quad\quad
 +21\left(\frac{1}{z_1}-\frac{1}{z_2}\right)
Q_{1,1}(z_2,z_2)=0, 
\nn
\lefteqn{\left(\frac{3}{7z_1^{2}}+\frac{3}{7z_2^{2}}
-21(\partial_{z_1}
+\partial_{z_2})^2-29
\left(\frac{\partial_{z_1}}{z_1}
+\frac{\partial_{z_2}}{z_2}\right)
\right)Q_{1,1}(z_1,z_2)}
\nn
&&\quad \quad +\left(\left(\frac{1}{z_2}-\frac{1}{z_1}\right)
\partial_{z_1}\partial_{z_2}
+\frac{1}{7}\left(\frac{\partial_{z_1}}{z_2^{2}}-
\frac{\partial_{z_2}}{z_1^{2}}\right)\right)Q_{0,0}(z_1,z_2)=0.
\end{eqnarray*}
Note that this is a system of regular singular points.

Since 
\begin{eqnarray*}
Q_{0,0}(z_1,z_2)&=&z_1^{1/14}z_2^{1/14}R_{0,0}((z_2/z_1)^{1/2}),\\
Q_{1,1}(z_1,z_2)&=&z_1^{1/14-1/2}z_2^{1/14-1/2}
R_{1,1}((z_2/z_1)^{1/2}),
\end{eqnarray*}
we obtain the following system
for $R_{0,0}(z)$ and $R_{1,1}(z)$:
\bea 
\lefteqn{R_{0,0}''(z)-\frac{-43z^{-1}+100z-54z^{3}}{21(1-z^2)^2}R_{0,0}'(z)}\nn
&&-\frac{23z^2+23z^{-2}-6}{147(1-z^2)^2}R_{0,0}(z)+\frac{4}{z(1-z^2)}R_{1,1}(z)=0,\label{eq3}\\
\lefteqn{R_{1,1}''(z)-\frac{100z-117z^{3}-z^{-1}}{21(1-z^2)^2}R_{1,1}'(z)}\nn
&& -\frac{324}{21(1-z^2)^2} R_{1,1}(z)+\frac{83z-203z^{3}+114z^{5}-21+42z^2-21z^4}{21z(z^2-1)^3}
R_{0,0}'(z) \nn
&&-\frac{37z^4+6z^2+37}{147z(z^2-1)^3}R_{0,0}(z)=0.\label{eq3-1}
\eea
Note that this system is analytic in the region $0<|z|<1$ and 
has a regular singularity at $z=0$.
Using the standard method, one can prove the convergence of power 
series solutions of (\ref{eq3})--(\ref{eq3-1}).

Similarly, from
\begin{eqnarray*}
\lefteqn{\langle w'_{(4)}, \mathcal{ Y}_{1}(G(-1/2)w_{(1)}, z_{1})
\mathcal{ Y}_{2}(w_{(2)}, z_{2})(8L(-2)-21L(-1)^2} \nn
&& \quad\quad\quad\quad\quad\quad\quad\quad\quad \quad\quad\quad
+21G(-3/2)G(-1/2))w_{(3)} \rangle=0, \nn
\lefteqn{\langle w'_{(4)}, \mathcal{ Y}_{1}(w_{(1)}, z_{1})
\mathcal{ Y}_{2}(G(-1/2)w_{(2)}, z_{2})(8L(-2)-21L(-1)^2} \nn
&& \quad\quad\quad\quad\quad\quad\quad\quad\quad \quad\quad\quad
+21G(-3/2)G(-1/2))w_{(3)} \rangle=0,
\end{eqnarray*}
we obtain the following system:
\begin{eqnarray*}
\lefteqn{\left(\left(\frac{3}{7z_1^{2}}
-\frac{4}{7z_2^{2}}\right)-
21(\partial_{z_1}+\partial_{z_2})^2-29
\left(\frac{\partial_{z_1}}{z_1}
+\frac{\partial_{z_2}}{z_2}\right)\right)Q_{1,0}(z_1,z_2)} \nn
&&\quad\quad\quad\quad\quad\quad \quad\quad\quad
 +\left(21\left(\frac{\partial_{z_1}}{z_2}-\frac{\partial_{z_1}}{z_1}\right)-
\frac{21}{7z_1^2}\right)
Q_{0,1}(z_2,z_2)=0, 
\nn
\lefteqn{\left(\left(\frac{3}{7z_2^{2}}
-\frac{4}{7z_1^{2}}\right)-
21(\partial_{z_1}+\partial_{z_2})^2-29
\left(\frac{\partial_{z_1}}{z_1}
+\frac{\partial_{z_2}}{z_2}\right)\right)Q_{0,1}(z_1,z_2)} \nn
&&\quad\quad\quad\quad\quad\quad \quad\quad\quad
 +\left(21\left(\frac{\partial_{z_2}}{z_1}-\frac{\partial_{z_2}}
{z_2}\right)-
\frac{21}{7z_2^2}\right)
Q_{1,0}(z_2,z_2)=0.
\end{eqnarray*}

Since
\begin{eqnarray*}
Q_{1,0}(z_1,z_2)&=&z_1^{1/14-1/2}z_2^{1/14}R_{1,0}((z_2/z_1)^{1/2}),\\
Q_{0,1}(z_1,z_2)&=&z_1^{1/14}z_2^{1/14-1/2}
R_{0,1}((z_2/z_1)^{1/2}),
\end{eqnarray*}
we obtain a system
for $R_{1,0}(z)$ and $R_{0,1}(z)$:
\bea
\lefteqn{R_{1,0}''(z)-\frac{168z
-118z^{3}-116z^{-1}}{42(z^{2}-1)^2}R_{1,0}(z)} \nn
&& -\frac{\frac{36}{147}-\frac{348}{147}z^{-2}+\frac{14700}{28182}}
{(z^{2}-1)^2}R_{1,0}(z)
-\frac{2}{z^{2}-1}R_{0,1}'(z)- \nn
&& \frac{2z^{-1}-6z}{7(z^{2}-1)^2}R_{0,1}(z)=0, \nn
\lefteqn{R_{0,1}''(z)+\frac{-425z^{-1}
+424z^{3}}{7(z^{2}-1)^2}R_{0,1}(z)} \nn
&& -\frac{12z^{-2}-16z^{2}}{147(z^{2}-1)^2}R_{0,1}(z)
-\frac{2}{z^{2}-1}R_{1,0}(z)+ \nn
&& \frac{-84z+96z^{-1}}{21(z^{2}-1)^2}R_{1,0}(z)=0. 
\eea

\noindent {\small \sc Institut des Hautes \'{E}tudes Scientifiques, 
Le Bois-Marie, 35, Route De Chartres, F-91440 Bures-sur-Yvette, 
France}

\noindent {\it and}

\noindent {\small \sc Department of Mathematics, Rutgers University,
110 Frelinghuysen Rd., Piscataway, NJ 08854-8019 (permanent address)}

\noindent {\em E-mail address}: yzhuang@math.rutgers.edu

\vskip 1em

\noindent {\small \sc Department of Mathematics, Rutgers University,
110 Frelinghuysen Rd., Piscataway, NJ 08854-8019}

\noindent {\em E-mail address}: amilas@math.rutgers.edu


\begin{thebibliography}{KWak2}

\bibitem[Ad]{A} 
D. Adamovi\'c, Rationality of Neveu-Schwarz vertex
operator superalgebras,  {\em Internat. Math. Res. Notices}
{\bf 1997} (1997), 865--874.

\bibitem[AM]{AM} 
D. Adamovi\'c and A.  Milas, Vertex
operator algebras associated 
to modular invariant representations for
$A^{(1)}_1$, {\em Math. Res. Lett.} {\bf 2} (1995),  563--575.

\bibitem[As]{As} 
A. Astashkevich, On the structure of Verma modules over
Virasoro and Neveu-Schwarz algebras, {\em Comm. Math. Phys.} 
{\bf 186}  (1997), 531--562.





\bibitem[Ba1]{B1} 
K. Barron,  A supergeometric interpretation of
vertex operator superalgebras,  {\em Internat. Math. Res. Notices}
{\bf 1996} (1996), 409--430.

\bibitem[Ba2]{B2}
K. Barron, The supergeometric interpretation of vertex
operator superalgebras, Ph.D. thesis, Rutgers University, 1996.

\bibitem[Bo]{B}
R.~E.~Borcherds,
Vertex algebras, Kac-Moody algebras, and the Monster,
{\em Proc. Natl. Acad. Sci. USA} {\bf 83} (1986), 3068--3071.

\bibitem[C1]{C1}
J. L. Cardy, Effect of boundary conditions on the operator 
content of two-dimensional
conformally invariant theories,
{\em Nucl. Phys.} {\bf B275} (1986),  200--218.

\bibitem[C2]{C2}
J. L. Cardy,   Boundary conditions, fusion rules and 
the Verlinde formula, {\em Nucl. Phys.}
{\bf B324} (1989),  581--596.

\bibitem[D]{D} 
L. Dixon, Some world-sheet properties of superstring
compactifications, on orbifolds and otherwise, in: {\em Superstrings,
Unified Theories and Cosmology 1987 (Trieste, 1987)}, ICTP
Ser. Theoret. Phys., Vol. 4, World Sci. Publishing, Teaneck, NJ, 1988,
67--126.

\bibitem[DGH]{DGH}
L. Dixon, P. Ginsparg and J. Harvey, 
Beauty and the beast: superconformal symmetry in a
Monster module,  {\em Comm. Math. Phys.} {\bf 119}
 (1988), 221--241

\bibitem[DGM]{DGM}
L. Dolan, P. Goddard and P. Montague, Conformal field theory of twisted 
vertex operators, {\em Nucl. Phys.} {\bf B338} (1990),
529--601.

\bibitem[DMZ]{DMZ}
C.  Dong, G.  Mason, Y. Zhu, Discrete series of the Virasoro 
algebra and the moonshine
module, in: {\em 
Algebraic groups and their generalizations: quantum and
infinite-dimensional methods (University Park, PA, 1991)}, 
Proc. Sympos. Pure Math., Vol. 56, 
Part 2, Amer. Math. Soc., Providence,
RI, 1994, 295--316.


\bibitem[FHL]{FHL}
I.~B. Frenkel, Y.-Z. Huang and J.~Lepowsky,
On axiomatic approaches to vertex operator algebras and modules,
preprint, 1989;
{\em Memoirs Amer. Math. Soc.} {\bf 104}, 1993.

\bibitem[FLM1]{FLM1}
I. B. Frenkel, J. Lepowsky and A. Meurman, A natural representation of
the Fischer-Griess monster with the modular function $J$ as character,
{\em Proc. Natl. Acad. Sci. USA} {\bf 81} (1984), 3256--3260.



\bibitem[FLM2]{FLM2}
I.~B. Frenkel, J.~Lepowsky, and A.~Meurman,
{\em Vertex operator algebras and the Monster},
Pure and Appl. Math., {\bf 134}, Academic Press, New York, 1988.

\bibitem[FZ]{FZ} 
I. B. Frenkel and Y.  Zhu, Vertex operator
algebras associated to representations of affine and Virasoro
algebras, {\em Duke Math. J.} {\bf 66} (1992),  123--168.

\bibitem[FS]{FS}
J. Fuchs and C. Schweigert, D-brane conformal field theory, 
{\tt hep-th/9801190}, to appear.


\bibitem[Ge1]{G1}
D. Gepner, 
Exactly solvable string compactifications on manifolds of 
$\mbox{rm SU}(N)$ holonomy,
{\em Phys. Lett.} {\bf  B199} (1987), 380--388. 

\bibitem[Ge2]{G2}
D. Gepner, 
Space-time supersymmetry in compactified string theory and 
superconformal models, 
{\em Nucl. Phys.} {\bf B296} (1988), 757--778. 

\bibitem[GKO]{GKO} 
P. Goddard, A. Kent and D. Olive, 
Unitary representations of the Virasoro and super-Virasoro
algebras, {\em Comm. Math. Phys.} {\bf 103} (1986), 105--119.

\bibitem[Gr]{Gr}
B. R. Greene, Aspects of quantum geometry, in:
{\em Mirror Symmetry III, Proceedings of the Conference on Complex 
Geometry and Mirror Symmetry, Montr\'{e}al, 1995,} Stud. Adv. Math.,
Vol. 10, Amer. Math. Soc., Internat. Press and Centre de Recherches 
Math., 1--67. 

\bibitem[GP]{GP}
B. R. Greene and M. R. Plesser, Duality in Calabi-Yau moduli space,
{\em Nucl. Phys.} {\bf B338} (1990), 15--37.

\bibitem[H1]{H1}
Y.-Z. Huang, A theory of tensor products for module categories
for a vertex operator algebra, IV, {\em J. Pure Appl. Alg.}, {\bf 100}
(1995), 173-216.

\bibitem[H2]{H2}
Y.-Z. Huang, Virasoro vertex operator algebras,
(nonmeromorphic) operator product expansion and the tensor product
theory, {\em J. Alg.} {\bf 182} (1996), 201--234.

\bibitem[H3]{H2.1}
Y.-Z. Huang, A nonmeromorphic extension of the 
Moonshine module vertex operator algebra, in: 
{\em Moonshine, the Monster, and related topics 
(South Hadley, MA, 1994)},  Contemp. Math., {\bf 193}, Amer. Math.
Soc., Providence, RI, 1996, 123--148.

\bibitem[H4]{H2.5}
Y.-Z. Huang, Intertwining operator algebras, genus-zero modular
functors and genus-zero conformal field theories,  in: {\em Operads:
Proceedings of Renaissance Conferences}, ed. J.-L. Loday,
J. Stasheff, and A. A. Voronov, Contemporary Math., {\bf  202},
Amer. Math. Soc., Providence, 1997,  335--355.

\bibitem[H5]{H3}
Y.-Z. Huang, {\em Two-dimensional conformal geometry and vertex operator
algebras}, Progress in Mathematics, Vol. 148,
Birkh\"{a}user, Boston, 1997.

\bibitem[H6]{H4}
Y.-Z. Huang, Genus-zero modular functors and intertwining operator
algebras, {\it Internat. J. Math.} {\bf 9} (1998), 845--863.

\bibitem[H7]{H6}
Y.-Z. Huang, Generalized rationality and a generalized Jacobi identity
for intertwining operator algebras,  to appear.

\bibitem[HL1]{HL1}
Y.-Z. Huang and J. Lepowsky,
\newblock Toward a theory of tensor products for representations of a 
vertex operator algebra, 
\newblock in: {\em Proc. 20th International Conference on Differential 
Geometric Methods in 
Theoretical Physics, New York, 1991}, 
\newblock ed. S. Catto and A. Rocha, 
\newblock World Scientific, Singapore, 1992, Vol. 1, 344--354.

\bibitem[HL2]{HL4}
Y.-Z. Huang and J. Lepowsky, Tensor products of modules for a vertex
operator algebra and vertex tensor categories, in:
     {\em Lie Theory and Geometry,
in honor of Bertram Kostant,}
ed. R. Brylinski, J.-L. Brylinski, V. Guillemin, V. Kac,
Birkh\"{a}user, Boston, 1994, 349--383.

\bibitem[HL3]{HL2}
Y.-Z. Huang and J. Lepowsky, A theory of tensor products for module
categories for a vertex operator algebra, I, {\em Selecta
Mathematica, New Series} {\bf 1} (1995), 699-756.

\bibitem[HL4]{HL3}
Y.-Z. Huang and J. Lepowsky, A theory of tensor products for module
categories for a vertex operator algebra, II, {\em Selecta
Mathematica, New Series} {\bf 1} (1995), 757--786.

\bibitem[HL5]{HL5}
Y.-Z. Huang and J. Lepowsky, A theory of tensor products for module
categories for a vertex operator algebra, III, {\em J. Pure
Appl. Alg.} {\bf 100} (1995),  141-171.

\bibitem[HL6]{HL5.5}
Y.-Z. Huang and J. Lepowsky, 
Intertwining operator algebras and vertex tensor categories for
affine Lie algebras, {\em Duke Math. J.}, to appear. 

\bibitem[HL7]{HL6}
Y.-Z. Huang and J. Lepowsky, A theory of tensor products for module
categories for a vertex operator algebra, V, 
to appear.





\bibitem[KWan]{KW} 
V. Kac and W. Wang,  Vertex operator 
superalgebras and their representations, 
in: {\em 
Mathematical aspects of conformal and topological field theories and 
quantum groups (South Hadley, MA,
1992)}, 
Contemp. Math., Vol. 175, 161--191.

\bibitem[KWak1]{KWak1} 
V. Kac and M. Wakimoto, Unitarizable highest
weight representations of the Virasoro, Neveu-Schwarz and Ramond
algebras, in: {\em Conformal groups and related symmetries: physical
results and mathematical background (Clausthal-Zellerfield, 1985),}
Lecture Notes in Phys. Vol. 261, Springer, Berlin, 1986, 345--371.

\bibitem[KWak2]{KWak2} 
V. Kac and  M. Wakimoto, Modular
invariant representations of infinite-dimensional Lie algebras and
superalgebras, {\em Proc. Nat. Acad. Sci. U.S.A.} {\bf 85}
(1988), 4956--4960.

\bibitem[LVW]{LVW}
W. Lerche, C. Vafa and N. P. Warner, Chiral rings in $N=2$ 
superconformal theories, {\em Nucl. Phys.} {\bf B324} (1989), 
427--474.

\bibitem[L]{L}
H. Li, Representation theory and tensor product theory for 
vertex operator algebras, PhD thesis, Rutgers University 1994.

\bibitem[RS]{RS}
A. Recknagel and V.  Schomerus, D-branes 
in Gepner models, {\em Nucl. Phys.} {\bf B531} (1998), 
185--225.

\bibitem[S1]{S1}
G. B.~Segal,
The definition of conformal field theory,
preprint, 1988.

\bibitem[S2]{S2}
G. B.~Segal, Two-dimensional conformal field theories and modular
functors, in: {\em Proceedings of the IXth International Congress on
Mathematical Physics, Swansea, 1988},
Hilger, Bristol, 1989, 22--37.


\bibitem[Z]{Z} 
Y. Zhu,  Modular invariance of characters of
vertex operator algebras, {\em J. Amer.  Math. Soc.}  {\bf 9}
(1996), 237--302.




\end{thebibliography}
\end{document}